\documentclass[11pt,reqno,a4paper]{amsart}
\usepackage{amsmath,amsthm,verbatim,amscd,amssymb,setspace,bbm,enumitem,hyperref}
\usepackage[small]{diagrams}
\usepackage{exscale,color,graphicx}
\usepackage[normalem]{ulem}
\usepackage{rotating,float}
\newcommand{\id}{\operatorname{id}}
\newcommand{\p}{\partial}
\newcommand{\oV}{\mkern 4mu\overline{\mkern-1mu V\mkern-1mu}\mkern 1mu}
\newcommand{\N}{\mathbb{N}}
\newcommand{\Z}{\mathbb{Z}}
\newcommand{\Q}{\mathbb{Q}}
\newcommand{\R}{\mathbb{R}}
\newcommand{\C}{\mathbb{C}}
\renewcommand{\P}{\mathbb{P}}

\renewcommand{\leq}{\leqslant}
\renewcommand{\geq}{\geqslant}
\renewcommand{\epsilon}{\varepsilon}

\newtheoremstyle{fancy}{}{}{\itshape}{}{\textbf\bgroup}{.\egroup}{ }{}
\newtheoremstyle{fancy2}{}{}{\rm}{}{\textbf\bgroup}{.\egroup}{ }{}

\theoremstyle{fancy}
\newtheorem{theorem}{Theorem}[section]
\newtheorem{lemma}[theorem]{Lemma}

\newtheorem{prop}[theorem]{Proposition}

\newcounter{mtheorem}
\newtheorem{mtheorem}[mtheorem]{Theorem}

\newtheorem{mcor}[mtheorem]{Corollary}

\theoremstyle{fancy2}
\newtheorem{definition}[theorem]{Definition}
\newtheorem{example}[theorem]{Example}
\newtheorem{remark}[theorem]{Remark}
\newtheorem{observation}[theorem]{Observation}

\numberwithin{equation}{section}

\begin{document}
\title{Asymptotically conical Calabi-Yau metrics \\ on quasi-projective varieties}
\date{\today}
\author{Ronan J.~Conlon}
\address{D\'epartement de Math\'ematiques, Universit\'e du Qu\'ebec \`a Montr\'eal, Case Postale 8888, Succursale
Centre-ville, Montr\'eal (Qu\'ebec), H3C 3P8, Canada}
\email{rconlon@cirget.ca}
\author{Hans-Joachim Hein}
\address{Department of Mathematics, University of Maryland, College Park, MD 20742--4015, USA}
\email{hein@umd.edu}
\date{\today}
\begin{abstract}
Let $X$ be a compact K\"ahler orbifold without $\C$-codimension-$1$ singularities. Let $D$ be a suborbifold divisor in $X$ such that $D \supset {\rm Sing}(X)$ and $-pK_X = q[D]$ for some $p, q \in \N$ with $q > p$. Assume that $D$ is Fano. We prove the following two main results. (1) If $D$ is K\"ahler-Einstein, then, applying results from our previous paper \cite{Conlon}, we show that each K\"ahler class on $X \setminus D$ contains a unique asymptotically conical Ricci-flat K\"ahler metric, converging to its tangent cone at infinity at a rate of $O(r^{-1-\epsilon})$ if $X$ is smooth. This provides a definitive version of a theorem of Tian and Yau \cite{Tian}. (2) We introduce new methods to prove an analogous statement (with rate $O(r^{-0.0128})$) when $X = {\rm Bl}_p \P^3$ and $D = {\rm Bl}_{p_1,p_2}\P^2$ is the strict transform of a smooth quadric  through $p$ in $\P^3$. Here $D$ is no longer K\"ahler-Einstein, but the normal $\mathbb{S}^1$-bundle to $D$ in $X$ admits an irregular Sasaki-Einstein structure which is compatible with its canonical CR structure. This provides the first example of an affine Calabi-Yau manifold of Euclidean volume growth with irregular tangent cone at infinity.
\end{abstract}
\maketitle
\markboth{Ronan J.~Conlon and Hans-Joachim Hein}{Asymptotically conical Calabi-Yau metrics on quasi-projective varieties}

\section{Introduction}

\subsection{An optimal Tian-Yau theorem}\label{s:intro:1} Our first main result optimises the following theorem due to Tian and Yau \cite[Corollary 1.1]{Tian}. Here we state it using slightly different terminology.\medskip

\noindent {\bf Tian-Yau Theorem.} \emph{Let $X^n$ be a compact K\"ahler orbifold without $\C$-codimension-$1$ singularities. Let $D \supset {\rm Sing}(X)$ be a neat and almost ample suborbifold divisor in $X$ with $-K_{X}=\alpha[D]$ for some $\alpha>1$. If $D$ admits K\"ahler-Einstein metrics of positive scalar curvature, then $X\setminus D$ admits complete Ricci-flat K\"ahler metrics $g$ of Euclidean volume growth.
Moreover, if we denote by $\rho$ the distance function from some fixed point with respect to $g$, then the curvature tensor of $g$ is $O(\rho^{-2})$ with respect to the $g$-norm; it is $o(\rho^{-2})$ if and only if $D$ is biholomorphic to $\P^{n-1}$.}\medskip

Before stating our improved version of this theorem, let us make some clarifying remarks.

$\bullet$ Proving such a result for orbifolds rather than only for manifolds requires little extra work but includes many more examples. Indeed, even to recover the well-known Ricci-flat ALE spaces in real dimension $4$, we require $X$ to be singular in all examples except three; compare \cite{rasdeaconu} and Appendix \ref{s:kronheimer}. These ALE spaces provide excellent examples for many typical orbifold phenomena.

$\bullet$ Assuming the absence of $\C$-codimension-$1$ singularities in $X$, i.e.~conical singularities along $D$, is no restriction of generality but greatly aids clarity in the statement and in the proof.

$\bullet$ The \emph{neat} and \emph{almost ample} conditions, introduced in \cite[Definition 1.1(i), (ii)]{Tian}, are complicated and technical, whereas the \emph{suborbifold} (or \emph{admissible}) condition naturally generalises the notion of a smooth complex hypersurface in a smooth complex manifold \cite[Definition 1.1(iii)]{Tian}.

$\bullet$ It is clear from \cite{Tian} that $\alpha$ is not meant to be restricted to be an integer. The condition that $-K_X = \alpha[D]$ thus requires some clarification due to the possibility of torsion in Pic$^{\rm orb}(X)$.

$\bullet$ The K\"ahler-Einstein metric on $D$ is understood to be an orbifold metric adapted to the orbifold structure on $D$ inherited from $X$. Thus, we may require the existence of K\"ahler-Einstein orbifold metrics with $\C$-codimension-$1$ singularities on $D$ even if $D$ is smooth as a variety.

$\bullet$ Asymptotically conical manifolds clearly have quadratic curvature decay and Euclidean volume growth. In the Ricci-flat case, there are difficult regularity theorems going in the opposite direction \cite{ALE, Cheeger, CM}. However, by analysing the proof in \cite{Tian}, one can prove directly that the Ricci-flat metrics obtained in \cite{Tian} are in fact asymptotically conical in the strongest possible sense.

$\bullet$ The final condition that $D \cong \P^{n-1}$ refers to an isomorphism of varieties, not of orbifolds.

Our first main result, to be proved in Section \ref{s:optimal_TY}, can now be summarised as follows.

\begin{mtheorem}\label{mainmain}
Let $X^n$ be a compact K\"ahler orbifold without $\C$-codimension-$1$ singularities. Let $D$ $\supset$ ${\rm Sing}(X)$
be a suborbifold divisor in $X$ such that $-pK_X = q[D]$ with $p,q \in \N$ and $\alpha = \frac{q}{p} > 1$, and such that $D$ admits a K\"ahler-Einstein metric of positive scalar curvature. For every K\"ahler class $\mathfrak{k}$ on $X \setminus D$ and for every $c>0$, there exists a unique
Ricci-flat K\"ahler metric $\omega_c \in \mathfrak{k}$ satisfying
\begin{equation}\label{e:main_rate_0}
|\nabla_{g_0}^k(\exp^{*}(g_{c})-cg_{0})|_{g_0 } \leq C(k)r^{-\lambda-k}
\end{equation}
for some $\lambda > 0$ and all $k \in \N_0$. Here $\exp:N_{D}\to X$ denotes the exponential map of any background Hermitian metric on $X$,
$g_{0}$ is the pullback of the Calabi ansatz Ricci-flat K\"ahler cone metric under the covering map $N_D \setminus 0 \to K_D^p\setminus 0$ induced by the adjunction isomorphism $N_D^{q-p} \cong K_D^{-p}$, and $r$ is the radius function of the cone metric $g_0$. Moreover, $\omega_c$ in fact satisfies \eqref{e:main_rate_0} with
\begin{equation}\label{e:main_rate}
\lambda = \min\{2-\epsilon, \frac{n}{\alpha-1}\}
\end{equation}
for any $\epsilon > 0$. As a direct consequence of the uniqueness, $\omega_c$ is invariant under all automorphisms of $(X,D)$ that preserve $\mathfrak{k}$ and induce isometries of $g_0$, and $\omega_{t c, t \mathfrak{k}} = t \omega_{c,\mathfrak{k}}$ for all $t > 0$.
\end{mtheorem}

Let us again make a few remarks to clarify the statement of this theorem.

$\bullet$ It may not be possible to cancel any common factors of $p$ and $q$ in the equation $-pK_X = q[D]$, even if $\alpha = \frac{q}{p} \in \N$; see Appendix \ref{s:cyclic_quotient}. However, $\pi_1(X) = 0$ from \cite{takayama} and Proposition \ref{t:vanishing}, so if $X$ is {smooth}, then ${\rm Pic}(X)$ is torsion-free and we can indeed assume that $p,q$ are coprime. On the other hand, if $X$ is singular, then $\pi_1^{\rm orb}(X)$ may be nontrivial, hence ${\rm Pic}^{\rm orb}(X)$ may contain torsion.\footnote{Cristiano Spotti pointed out to us that this phenomenon was first observed in \cite{Mumford} for the classical elliptic modular curve $\mathbb{H}/{\rm PSL}(2,\Z)$, which is isomorphic to $\C$ as a variety but whose $\pi_1^{\rm orb}$ and ${\rm Pic}^{\rm orb}$ are nontrivial.}

$\bullet$ For us, a \emph{K\"ahler class} on an open complex manifold is simply a de Rham cohomology class (of degree $2$ and with real coefficients) that contains closed positive $(1,1)$-forms.

$\bullet$ Note that $\omega_c$ is only claimed to be unique among K\"ahler forms in $\mathfrak{k}$ that satisfy (\ref{e:main_rate_0}) for some $\lambda > 0$. This uniqueness is an immediate consequence of \cite[Theorem 3.1]{Conlon}.

$\bullet$ Recall that the \emph{Calabi ansatz} produces a Ricci-flat K\"ahler cone metric on $K_D \setminus 0$ (with apex at the zero section) from a given K\"ahler-Einstein metric on $D$; see \cite[Section 1.3.3]{Conlon}  and Section \ref{s:orbi_calabi}. This metric is U$(1)$-invariant and hence pushes down from $K_D \setminus 0$ to $K_D^p \setminus 0$.

$\bullet$ The K\"ahler-Einstein metric on $D$ is unique only up to the action of $G = {\rm Aut}(D)_0$. Thus, if $G$ is nontrivial, then we actually have a family of Ricci-flat metrics $\omega_{c,x} \in \mathfrak{k}$, parametrised by the points $x$ of the symmetric space associated with $G$ and differing from each other by $O(1)$ at infinity. It is not clear (and probably false in general) that these metrics only differ by automorphisms of $X \setminus D$, but they may still be isometric; compare \cite[Remark 2.6, Corollary 3.14, Remark 5.10]{Conlon}.\medskip

\noindent {\bf Comparing the two theorems and their proofs.} Many of the refinements in Theorem \ref{mainmain} (no neat or almost ample condition, all K\"ahler classes, the parameter $c$, uniqueness and symmetry) are due to an  improvement of general technique in \cite{Conlon}, partly based on important earlier contributions by van Coevering \cite{vanC}, whereas asymptotics of the form (\ref{e:main_rate_0}) are already implicit in Tian-Yau \cite{Tian}. Let us point out one useful consequence of our explicit estimate
(\ref{e:main_rate}).

\begin{mcor}\label{c:rate}
If $X$ is smooth, then the best possible convergence rate $\lambda$ of the Ricci-flat metrics of Theorem \ref{mainmain} to their tangent
cones at infinity is always strictly greater than $1$.
\end{mcor}

\begin{proof}
Since $N_D^{q-p} = K_D^{-p}$ and ${\rm Pic}(D)$ is torsion-free because $\pi_1(D) = 0$, there exists a line bundle $L$ with $L^p = N_D$ and $L^{q-p} = K_D^{-1}$. Thus, by \cite[p.~32, Corollary]{KO}, $q - p \leq n$ with equality if and only if $D = \P^{n-1}$, so that $\alpha - 1 \leq n$ and $\lambda \geq 1$ with equality if and only if $D = \P^{n-1}$, $N_D = \mathcal{O}(1)$. But in the latter case, the cone, and hence $(X \setminus D, g_c)$ itself, must be isometric to flat $\C^n$.

It seems reasonable to expect that $\alpha - 1 \leq n$ even if $X$ is singular. Moreover, equality should still imply that the cone is $\C^n/\Gamma$ (see Remark \ref{r:obata} for some examples where $\Gamma \neq \{1\}$), so that $\lambda \geq 2n$ by \cite[Theorem 5.103]{Cheeger}. In fact, if $N_D$ is effective, then the inequality follows from \cite[Theorem 8.7]{ross-thomas}. Alternatively, under the weaker assumption that $N_D$ has a multivalued section defining a covering space of $N_D^*\setminus 0$, the inequality as well as the equality case follow from \cite[Sections 2.2--2.3]{GMSY:obs}.
\end{proof}

We also mention a curious technical detail. One key point in both Tian-Yau's proof \cite{Tian} and ours is the existence of a nonnegatively curved Hermitian metric $h$ on the line orbibundle $[D]$ that has strictly positive curvature in a neighbourhood of $D$. Tian-Yau used the almost ampleness of $D$ to construct $h$, and then used $h$ to construct reference metrics for solving a Monge-Amp{\`e}re equation on $X \setminus D$. We construct $h$ by hand ({deducing} almost ampleness as a \emph{corollary} of the existence of $h$ by applying some deep results of Grauert, but not using this almost ampleness in itself) and then use $h$ to prove that $h^{i,0}(X) = 0$ for all $i > 0$, the case $i = 2$ being precisely what is required to apply the theory of \cite{Conlon}. This vanishing was also noted by Tian and Yau, but not in their existence
proof: in \cite[Section 6]{Tian}, following a classical idea of Kobayashi \cite{Kobayashi}, they used it to show that $\pi_1(X) = 0$. As mentioned above, $\pi_1(X) = 0$ was later proved much more generally in \cite{takayama}.

\subsection{An affine Calabi-Yau manifold with irregular tangent cone} If we have $-pK_X = q[D]$ with $\alpha = \frac{q}{p} > 1$, and if we only assume that $D$ is Fano (or equivalently, that $N_D$ is positive) rather than K\"ahler-Einstein of positive scalar curvature, then in general there is no obvious candidate for an asymptotically Ricci-flat model metric defined in a punctured tubular neighbourhood of $D$.

In \cite[Example 6.2]{vanC}, van Coevering pointed out the example $X = {\rm Bl}_p\P^3$, with $D \in |{-\frac{1}{2}K_X}|$ the proper transform of a smooth quadric in $\P^3$ passing through $p$. Then $D = {\rm Bl}_{p_1,p_2}\P^2$, which is Fano but does not admit any K\"ahler-Einstein metrics. However, it is known from \cite{futaki} that $N_D \setminus 0$ admits a toric Calabi-Yau cone metric $\omega_0 = \frac{i}{2}\partial\bar\partial r^2$, with infinite end at the zero section, whose Reeb vector field $J(r \partial_r)$ is irregular, i.e.~does not generate an $\mathbb{S}^1$-action. It is natural to ask whether a Tian-Yau theorem can be proved based on this irregular Calabi-Yau cone structure on $N_D \setminus 0$.

This turns out to be a difficult question because $X \setminus D$ is not a resolution of singularities of the cone but, on the contrary, an affine variety. Thus, it is not even clear that the complex structure of $X \setminus D$ converges to the complex structure of the cone with respect to $g_0$ and any reasonable choice of a diffeomorphism at infinity between the two spaces.  Experience with Sasakian geometry indeed suggests that it need not. However, in Section \ref{s:irregular}, we will prove the following.

\begin{mtheorem}\label{t:second_main}
Let $r$ denote the radius function of the cone metric $g_0$. There exists a diffeomorphism $\Phi$ from a neighbourhood of the zero section of $N_D$ onto a neighbourhood of $D$ in $X$
such that for all $\mathfrak{k} \in H^2(X \setminus D) \cong \R$ and all $c > 0$ there exists a Calabi-Yau metric $\omega_c \in \mathfrak{k}$ such that
\begin{equation}\label{e:second_main}
|\nabla^{k}_{g_0}(\Phi^{*}g_{c}-c g_{0})|_{g_{0}} \leq C(k) r^{-0.0128 - k}
\end{equation}
for all $k \in \N_0$. If $\mathfrak{k} = 0$, then we can improve $0.0128$ slightly to $0.0192$. Furthermore, the $\omega_c$ are the only Calabi-Yau metrics in $\mathfrak{k}$ that are asymptotically conical with respect to $\Phi$, so that $\omega_{t c, t \mathfrak{k}} = t \omega_{c,\mathfrak{k}}$ for all $t > 0$, and they are invariant under the maximal compact $(\mathbb{S}^1)^2 \subset {\rm Aut}(X,D) = (\C^*)^2$.
\end{mtheorem}

This provides the first example of a Calabi-Yau space of Euclidean volume growth with irregular asymptotic cone which, as a complex manifold, is not a crepant resolution of its asymptotic cone. Currently, we are not aware of any other candidates. See Section \ref{s:vanCapp} for some more discussion.

Regarding the statement of the theorem, let us quickly note the following.

$\bullet$ The map $\Phi$ is not explicit. In fact, we will first prove a somewhat weaker result (Theorem \ref{irr-exi}), with a better rate for $k = 0$ but suboptimal derivative behaviour, where the map from the cone to the manifold is explicit; given this, we then construct $\Phi$ by an abstract gauge fixing.

$\bullet$ We have the same caveats regarding uniqueness as in Theorem \ref{mainmain}. In particular, \emph{asymptotically conical} is to be understood in the strict sense of \cite[Definition 1.11]{Conlon}, and the same existence result holds with $\omega_0$ replaced by any of its pullbacks under ${\rm Aut}(D) = (\C^*)^2$ acting on the cone.

Finally, we mention that, whereas the proof of Theorem \ref{irr-exi} is similar to the proof of Theorem \ref{mainmain} in broad outline, there are many additional technical difficulties stemming from the irregularity of the cone metric and the fact that $X \setminus D$ is not biholomorphic to the cone away from a compact set.\footnote{The key issue is that smooth functions on $X$ are \emph{not polyhomogeneous} with respect to the irregular radius function, no matter what diffeomorphism we use to identify $X \setminus D$ with the cone at infinity; cf. \cite[Section 5]{Rochon}.} We solve these problems, relying on a detailed understanding of $(X,D)$ and $g_0$ and some fortuitous numerology:~if the {weights (\ref{irr:weights}) of the irregular Reeb field $J(r\partial_r)$ had been different by as little as $1\%$, then we would have been unable to prove Theorem \ref{t:second_main} for several reasons.}

\subsection{Notation and terminology} Most of our notation is fairly standard except possibly for $L \setminus 0$ for the total space of a line bundle $L$ minus its zero section, and $L^\times$ for a  negative holomorphic line bundle with its zero section contracted to a point. We refer to \cite[Section 1.3]{Conlon} for some specialised terminology and background used throughout this project.
Our main resource for orbifold geometry is \cite[Chapter 4]{book:Boyer} and our only definition not found in this chapter is that of a \emph{suborbifold divisor}:~by this, we mean an \emph{admissible divisor} in the sense of Tian and Yau; in particular, it is the $\Q$-Cartier divisor associated with a suborbifold of complex codimension $1$.
 Also, we consistently work within the orbifold category:~all line bundles are really line orbibundles, and $K_Z$ really means $K_Z^{\rm orb}$.

\subsection{Acknowledgements}
This paper is a continuation of our previous article \cite{Conlon}, and we would like to thank the same people and institutions as in \cite{Conlon}. In particular, Theorem \ref{mainmain} is directly based on parts of RJC's thesis (Imperial College 2011) supervised by Mark Haskins, and HJH thanks the EPSRC for postdoctoral support under Leadership Fellowship EP/G007241/1 until 2013. RJC also wishes to acknowledge support from a Britton postdoctoral fellowship held at McMaster University until 2013, and he thanks Tristan Collins for explaining the method of \cite{Tristan} as well as Kael Dixon for
providing Figure \ref{flow-pic}. HJH also thanks the CIRGET/UQAM for hosting him in September 2013 and Song Sun for a useful conversation related to \cite{ross-thomas}. \newpage
\section{Optimal existence if $D$ is K\"ahler-Einstein}\label{s:optimal_TY}

In this section, we provide the proof of Theorem \ref{mainmain} modulo the theory of \cite{Conlon}. We also omit some technical constructions; these we defer to Appendices \ref{app:c1}--\ref{s:gysin}. Assume for now that $p = 1$, so that $X \setminus D$ admits a holomorphic volume form $\Omega$ with a pole of order $\alpha = q > 1$ along $D$.

The first step is to endow $N_D\setminus 0$ with a natural Calabi-Yau cone structure $(g_0, \Omega_0)$, with infinite end at the zero section $0 \subset N_D$, using the given K\"ahler-Einstein metric on $D$; see Section \ref{s:covering}.

The second step is to show that $\Omega$ converges to $\Omega_0$ at a definite polynomial rate with respect to the cone metric $g_0$ and a suitable exponential map that we use to identify a tubular neighbourhood of the zero section in $N_D$ with a tubular neighbourhood of $D$ in $X$. See Section \ref{s:asy_hol_vf} for this.

We then complete the proof by appealing to our existence result from \cite{Conlon}. However, in order to apply our result, we must show that all K\"ahler classes on $X\setminus D$ are ``almost compactly supported''. This we prove in Section \ref{complete}, relying on the fact that $h^{2,0}(X)=0$ proved in Section \ref{hippers}.

Finally, in Section \ref{s:torsion_case}, we explain how to generalise our results from $p = 1$ to $p \in \N$.

\subsection{A Calabi-Yau cone structure on $N_D\setminus 0$.}\label{s:covering} Using the adjunction formula, we first construct a covering map from $N_D \setminus 0$ onto $K_D \setminus 0$. The orbifold Calabi ansatz produces a Calabi-Yau cone structure on $K_D\setminus 0$, which we then pull back under this covering map.

\subsubsection{Construction of a covering map}\label{ss:covering} Let $(\tilde{U}, \Gamma, \varphi)$ be a uniformising chart for $X$ such that $\varphi(\tilde{U})$ meets $D$. Set
$\tilde{D} =\varphi^{-1}(D\cap\varphi(\tilde{U}))$. Let $\tilde{s}$ be a holomorphic defining function for $\tilde{D}$ on $\tilde{U}$ and let $\p_{\tilde{s}}$ denote the unique local trivialising section of $N_{\tilde{D}}$ on which $d\tilde{s}|_{\tilde{D}}$ evaluates to one. Define $\tilde{\Omega} = \varphi^*\Omega$ and observe that this has the following two properties.

$\bullet$ $\tilde\Omega$ is a $\Gamma$-invariant holomorphic volume form on $\tilde{U}\setminus\tilde{D}$ that blows up to order $\alpha$ along $\tilde{D}$.

$\bullet$ If $\lambda:(\tilde{U},\Gamma,\varphi)\to
(\tilde{U}',\Gamma',\varphi')$ is an injection of uniformising charts of $X$, then $\lambda^{*}\tilde\Omega'
=\tilde{\Omega}$ off $\tilde{D}$.

\noindent We can therefore construct a covering map $\tilde{\eta} : N_{\tilde{D}}\setminus 0 \to K_{\tilde{D}} \setminus 0$ by sending, for all $t \in \C^*$,
\begin{equation}\label{popcornless}
t\p_{\tilde{s}} \mapsto
t^{1-\alpha}\partial_{\tilde{s}}\,\llcorner\, (\tilde{s}^{\alpha}\tilde\Omega)|_{\tilde{D}}.
\end{equation}

Now the following easy observation is key: $\tilde{\eta}$ does not depend on our choice of defining function for $\tilde{D}$. Indeed, this immediately implies
that $\tilde{\eta}$ is $\Gamma$-equivariant (since $\tilde{s} \circ \gamma$ is a defining function for every $\gamma \in \Gamma$) and patches up across uniformising charts (since $\tilde{s}' \circ \lambda$ is a defining function for every injection $\lambda$ as above). Thus, we obtain a global covering map $\eta: N_D \setminus 0 \to K_D \setminus 0$, as desired.

\subsubsection{Orbifold Calabi ansatz}\label{s:orbi_calabi} The space $K_D \setminus 0$ carries a canonical holomorphic volume form $\Omega_\sharp$,
defined as follows. Let $(\tilde{U},\Gamma,\varphi)$ be a local uniformising chart for $D$ with coordinates $(v_1, \dots, v_{n-1})$. Define local coordinates on $K_{\tilde{U}}$ by $(v_1,\dots, v_n) \mapsto v_n dv_1 \wedge \ldots \wedge dv_{n-1}$.  Then the holomorphic volume form $\tilde{\Omega}_\sharp = dv_1 \wedge \ldots \wedge dv_n$ on $K_{\tilde{U}}$ is $\Gamma$-invariant and patches up across uniformising charts.

Now let $\omega$ be a K\"ahler-Einstein metric on $D$ with ${\rm Ric}(\omega) = \omega$. Then $\tilde{\omega} = \varphi^*\omega$ defines a Hermitian fibre metric $\tilde{h}$ on $K_{\tilde{U}}$ via $\tilde{h}(v, v) \tilde{\omega}^{n-1} = i^{(n-1)^2}v \wedge \bar{v}$. By the Calabi ansatz \cite[Proposition 3.1]{Lebrun}, the function $\tilde{s}$ given by $\tilde{s}^{2n} = \tilde{h}$ is the radius function of a Calabi-Yau cone metric $\tilde\omega_{\sharp}$ on  $K_{\tilde{U}} \setminus 0$ such that $\tilde{\Omega}_{\sharp}$ has constant norm with respect to $\tilde{\omega}_{\sharp}$. Thus, we obtain a Calabi-Yau cone structure $(\omega_{\sharp}, \Omega_{\sharp})$ on $K_D \setminus 0$, with apex at the zero section $0 \subset K_D$, such that $(2n)^n\omega_{\sharp}^n = i^{n^2}\Omega_{\sharp} \wedge \bar{\Omega}_{\sharp}$.

\subsection{Asymptotics of the holomorphic volume forms}\label{s:asy_hol_vf} Now put $\Omega_0 = (-1)^n(\alpha-1)^{-1}\eta^*\Omega_{\sharp}$ and $\omega_0 = 2n(\alpha-1)^{-2/n}\eta^*\omega_{\sharp}$. Then $(\omega_0, \Omega_0)$ is a Calabi-Yau cone structure on $N_D \setminus 0$, with infinite end at the zero section $0 \subset N_D$, satisfying the standard normalisation $\omega_0^n = i^{n^2}\Omega_0 \wedge \bar\Omega_0$.

\begin{prop}\label{rate}
Fix a background Hermitian metric $g$ on $X$, let $\exp:(T^{1,0}D)^\perp \to X$ denote the normal exponential map associated with $g$, and use orthogonal projection with respect to $g$ to identify $N_D$ and $(T^{1,0}D)^\perp$ as smooth complex line bundles. Then we have that $$\exp^{*}(\Omega)-\Omega_{0}=O(r^{-\frac{n}{\alpha-1}})\;\,\textrm{with $g_{0}$-derivatives}$$
in a neighbourhood of $0 \subset N_D$. Here, $r$ denotes the radius function of the cone metric $g_0$.
\end{prop}

For clarity, we break the proof down into several steps.

\subsubsection{Convenient coordinates}
Let $(\tilde{U},\Gamma,\varphi)$ be a uniformising chart for $X$ near $D$ with holomorphic coordinates $(z_{1},\ldots,z_{n})$ such that $\tilde{D} =\varphi^{-1}
(D\cap\varphi(\tilde{U}))=\{z_{n}=0\}$ and
\begin{equation}\label{holo2}
\tilde{\Omega} =\varphi^{*}\Omega = z_n^{-\alpha} dz_{1}\wedge\ldots\wedge dz_{n}.
\end{equation}
We have associated holomorphic coordinates $(v_1, \ldots, v_n)$ on the total space of $K_{\tilde{D}}$,
representing the form $v_n dz_1 \wedge \ldots \wedge dz_{n-1}$ at $(v_1, \ldots , v_{n-1},0) \in \tilde{D}$. Similarly, we take $(w_1, \ldots, w_n)$ to be holomorphic coordinates on the total space of $N_{\tilde{D}}$ representing the coset of normal vectors $w_n(\partial_{z_n} + T^{1,0}\tilde{D})$ at $(w_1, \ldots, w_{n-1}, 0) \in \tilde{D}$. The covering map $\tilde{\eta}$ of Section \ref{ss:covering} may now be written as
\begin{equation*}\label{pullpop1}
(v_{1},\ldots,v_{n})= (w_{1},\ldots,w_{n-1},(-1)^{n-1}w_{n}^{1-\alpha}).
\end{equation*}
Consequently, bearing in mind the fact that $\tilde\Omega_\sharp = dv_1 \wedge \ldots \wedge dv_n$, we obtain that
\begin{equation}\label{holo23}
\tilde{\Omega}_{0}= (-1)^n(\alpha-1)^{-1}\tilde{\eta}^*\tilde\Omega_{\sharp} = w_n^{-\alpha}
dw_{1}\wedge\ldots\wedge dw_{n}.
\end{equation}

\subsubsection{Pulling back by the exponential map}\label{s:reg:exp} Abusing notation, we will write $\exp$ for the composition $N_D \to (T^{1,0}D)^\perp \to X$ of the $g$-orthogonal projection with the $g$-normal exponential map. Then $\exp$ is an orbifold diffeomorphism onto its image in some open neighbourhood of the zero section of $N_D$. Shrinking $U$ and lifting $\exp$, we obtain a $\Gamma$-equivariant diffeomorphism $\widetilde{\exp}: \tilde{V} \to \tilde{U}$.

Formulas (\ref{holo2}) and (\ref{holo23}) already
suggest that $\widetilde{\exp}^*(\tilde{\Omega}) - \tilde{\Omega}_0$ will be of lower order, but we require the results of Appendix \ref{app:c1} to make this precise. Thus, we begin by applying Observation \ref{obs:app_c1_first} with $\Phi = \widetilde{\exp}$; see also Example \ref{ex:app_c1}. This provides us with a new, smooth (but rarely ever holomorphic) coordinate system $z_i' = z_i + O(|z_n|)$ on $\tilde{U}$ such that, from (\ref{holo2}) and (\ref{eq:app_vol_form}),
\begin{equation}\label{exp2}
\tilde{\Omega} = \frac{dz_1' \wedge \ldots \wedge dz_n'}{(z_n')^\alpha} + \frac{\Upsilon \wedge dz_n'}{(z_n')^{\alpha-1}}
\end{equation}
for some smooth $(n-1)$-form $\Upsilon$ on $\tilde{U}$. On the other hand, by
Lemma \ref{lem:app_c1},
\begin{equation}\label{exp}
\widetilde{\exp}^{*}z'_{i}=w_{i}+A_{i,1} w_{n}^{2}+ A_{i,2} w_{n}\bar{w}_n+ A_{i,3}\bar{w}_{n}^{2}
\end{equation}
with $A_{i,j}$ smooth on $\tilde{V}$. Our bound on $\widetilde{\exp}^*(\tilde{\Omega}) - \tilde{\Omega}_0$ will follow from (\ref{exp2}), (\ref{exp}), and (\ref{holo23}).

\subsubsection{Estimates for smooth functions in terms of the cone metric}\label{s:regular:scaling} In order to proceed, we need a method of estimating the norm of smooth functions and of their derivatives on $\tilde{V}$ in terms of the cone metric $\tilde{g}_0$ and its radius function $\tilde{r} = \sqrt{2n}(\alpha-1)^{-1/n} \tilde\eta^* \tilde{s}$. We argue here using scaling.

Let $L$ denote the link of the cone $(N_{D}\setminus 0,g_{0}) \cong (\R^+ \times L, dr^2 \oplus r^2 g_L)$. As in \cite[Section 1.3.1]{Conlon}, we introduce the scaling diffeomorphism $\nu_{t}:[1,2]\times L\to [t,2t]\times L$  defined by $(r,x) \mapsto (tr,x)$ for $t > 0$. Lifting this map to our local uniformising chart for $N_{D}$, there must exist some $\mu > 0$ such that
$$\tilde{\nu}_{t}(w_{1},\ldots,w_{n})=(w_{1},\ldots,w_{n-1}, \mu w_{n}).$$
Using Section \ref{s:covering} and ignoring irrelevant constant factors, we can determine $\mu$ as follows:
\begin{equation*}\label{niceness}
t\tilde{r}(w_1, \ldots, w_n) =\tilde{r}(w_{1},\ldots,w_{n-1},\mu w_{n})= \tilde{s}(w_{1},\ldots,w_{n-1},\mu^{1-\alpha}w_{n}^{1-\alpha}) =\mu^{\frac{1-\alpha}{n}}\tilde{r}(w_{1},\ldots,w_{n}).
\end{equation*}
It subsequently follows from \cite[Lemma 1.6]{Conlon} that we have
\begin{eqnarray}
\label{reg_est1}
w_i = O(1)\;\,\textrm{with $\tilde{g}_{0}$-derivatives for all $i < n$},\\
\label{reg_est2}
w_{n}^{\pm 1}=  O(\tilde{r}^{\mp \frac{n}{\alpha-1}})\;\,\textrm{with $\tilde{g}_{0}$-derivatives}.
\end{eqnarray}

Finally, we prove that \emph{every} smooth function $A$ on $\tilde{V}$ satisfies the bounds (\ref{reg_est1}). Indeed,
\begin{equation*}
dA= \frac{\p A}{\p w_{n}}dw_{n}+\frac{\p A}{\p\bar{w}_{n}}d\bar{w}_{n} + \sum_{i=1}^{n-1}\frac{\p A}{\p w_i}dw_i +\sum_{i=1}^{n-1}\frac{\p A}{\p \bar{w}_i}d\bar{w}_i,
\end{equation*}
which, using (\ref{reg_est1}) and (\ref{reg_est2}), shows that $|dA|_{\tilde{g}_{0}}=O(\tilde{r}^{-1})$. More generally, by induction,
\begin{equation}\label{elephant}
A=O(1)\;\,\textrm{with $\tilde{g}_{0}$-derivatives}.
\end{equation}

\subsubsection{Asymptotics of the holomorphic volume forms}\label{s:reg:asympt} First of all, (\ref{exp2}) tells us that
\begin{equation*}
\widetilde{\exp}^{*}(\tilde{\Omega})=\frac{\widetilde{\exp}^{*}(dz_{1}')\wedge\ldots\wedge\widetilde{\exp}^{*}(dz_{n}')}{(\widetilde{\exp}^{*}z_{n}')^{\alpha}} + \frac{\widetilde{\exp}^{*}(\Upsilon) \wedge \widetilde{\exp}^{*}(dz_n')}{(\widetilde{\exp}^{*}z_n')^{\alpha-1}} = ({\rm I}) + ({\rm II}).
\end{equation*}
Using (\ref{exp}), (\ref{elephant}), (\ref{reg_est2}), and (\ref{holo23}), we see that
\begin{equation*}
\begin{split}
({\rm I}) &=[(dw_{1}+O(\tilde{r}^{-\frac{2n}{\alpha-1}-1}))\wedge\ldots\wedge (dw_{n}+O(\tilde{r}^{-\frac{2n}{\alpha-1}-1}))]/[w^{\alpha}_{n}(1+
O(\tilde{r}^{-\frac{n}{\alpha-1}}))^{\alpha}]\\
&=(1+O(\tilde{r}^{-\frac{n}{\alpha-1}}))(\tilde\Omega_{0}+ w_n^{-\alpha}O(\tilde{r}^{-\frac{2n}{\alpha-1}-1})O(\tilde{r}^{-1})^{n-1}) \\&=\tilde{\Omega}_{0}+O(\tilde{r}^{-\frac{n}{\alpha-1}}).
\end{split}
\end{equation*}
To deal with (II), we observe that $\widetilde{\exp}^*(\Upsilon)$ is a linear combination of wedge products of $n-1$ of the $2n$ basic forms $dw_1, d\bar{w}_1, \ldots, dw_n, d\bar{w}_n$, the coefficients being smooth functions on $\tilde{V}$. Thus, (\ref{elephant}), (\ref{reg_est1}), and (\ref{reg_est2}) show that $\widetilde{\exp}^*(\Upsilon) = O(\tilde{r}^{-1})^{n-1}$. Hence, using (\ref{exp}), (\ref{elephant}), and (\ref{reg_est2}),
$$({\rm II}) = O(\tilde{r}^{-1})^{n-1}O(\tilde{r}^{-\frac{n}{\alpha-1}-1}) O(\tilde{r}^{-\frac{n}{\alpha-1}})^{1-\alpha} = O(\tilde{r}^{-\frac{n}{\alpha-1}}).$$
This clearly finishes the proof of Proposition \ref{rate}.

\subsection{No holomorphic forms on $X$}\label{hippers}

Applying Proposition \ref{rate} and \cite[Theorem 2.4, Remark 2.5, and Theorem 3.1]{Conlon}, we find that Theorem \ref{mainmain} holds for every K\"ahler class on $X \setminus D$ that is $\mu$-almost compactly supported in the sense of \cite[Definition 2.3]{Conlon} for some $\mu < 0$. Precisely, if $n = 2$, then we also require that $K = \emptyset$ in this definition; moreover, in every dimension, we require that $\mu \leq -2$ in order to obtain the specific rate claimed in (\ref{e:main_rate}). In Section \ref{complete}, we will see that \emph{all} K\"ahler classes on $X\setminus D$ satisfy all of these conditions, thus completing the proof of Theorem \ref{mainmain} for $p = 1$.

The main ingredient for Section \ref{complete} is the following vanishing theorem.
\begin{prop}\label{t:vanishing}
Let $X$ be a compact K\"ahler orbifold with ${\rm Sing}(X)$ of complex codimension $> 1$.
Let $D$ be a suborbifold divisor in $X$ containing ${\rm Sing}(X)$ and satisfying $-K_{X}=q[D]$ for some $q\in\mathbb{N}$. Furthermore,
assume that the normal orbibundle to $D$ is positive. Then $h^{i,0}(X)=0$ for all $i>0$. In particular, it follows that $X$ is projective algebraic.
\end{prop}

The key point is to prove that the line orbibundle $[D]$ admits a nonnegatively curved Hermitian metric with strictly positive curvature on some tubular neighbourhood of $D$; see Lemma \ref{heythere} below.
Given this, an elementary generalisation of the Kodaira vanishing theorem due to Riemenschneider \cite[Theorem 6]{Riemen2}, whose proof works verbatim for K\"ahler orbifolds, {tells us that
$$h^{i,0}(X) =  h^{0,i}(X) =h^{i}(X,\mathcal{O}_X) = h^{i}(X,K_X + q[D]) = 0$$
for all} $i > 0$.\footnote{If $X$ is smooth and projective, then it follows directly from our assumptions that $[D]$ is big and nef, so we obtain the same conclusion from the Kawamata-Viehweg vanishing theorem for projective manifolds \cite[Section 6.D]{DemPCMI}.} The case $i = 2$ yields that $X$ is projective by Kodaira-Baily embedding \cite{Baily}.

\begin{lemma}\label{heythere}
The line orbibundle $[D]$ in Proposition \ref{t:vanishing} admits a nonnegatively curved Hermitian metric with strictly positive curvature on some tubular neighbourhood of $D$.
\end{lemma}

\begin{proof}
For a smooth divisor with positive normal bundle in a compact complex manifold, it is shown in \cite[\S VIII.1]{Griffiths} that it is always possible to find a tubular neighbourhood where the associated line bundle admits a metric of positive curvature. The proof in \cite{Griffiths} works for orbifolds as well. Thus, let $h$ be a Hermitian metric on $[D]$ with positive curvature on some open neighbourhood $U$ of $D$.

Fix once and for all a defining section $s$ of $[D]$. Let $f: \R^+ \to \R^+$ be a smooth function such that $f(t) = t(1 + g(t))$, where $g$ extends smoothly to $[0,\infty)$ with $g(0) = 0$. Then we can construct a new Hermitian metric $h_f$ on $[D]$ by setting $h_f(s,s) = f(h(s,s))$. The curvature of $h_f$ is given by
\begin{equation*}\label{janice}
-i\partial\bar{\partial}\log h_f=\biggl(\biggl(\frac{f'}{f}\biggr)^{2}-\frac{f''}
{f}
-\frac{f'}{hf}\biggr)i\partial h\wedge\bar{\partial}h+\frac{hf'}{f}
(-i\p\bar{\p}\log h),
\end{equation*}
where all Hermitian metrics are to be evaluated at $(s,s)$ and $-i\partial\bar\partial\log h > 0$ on $U$.

Fix $a,\delta \in \R^+$ such that the connected component containing $D$ of the open set $\{h(s,s) < 2a\delta\}$ is contained in $U$, and define $f(t) =\delta F(\frac{t}{\delta})$, where $F(x) = 1$ for $x \in [a,\infty)$ and $F(x) = G(H(x))$ for $x \in (0,a)$ with $G(H) = \frac{H}{1+H}$ and $H(x) = x \exp(\frac{x}{a-x})$. It remains to prove that
\begin{equation*}\label{e:ODI}
\biggl(\frac{F'(x)}{F(x)}\biggr)^{2}-\frac{F''(x)}{F(x)}
-\frac{F'(x)}{xF(x)}\geq 0\;\,{\rm for\;all}\;x \in (0,a).
\end{equation*}
In Appendix \ref{s:test_function}, we will see that this holds for all $a \geq 6$.
\end{proof}

The argument used to prove Proposition \ref{t:vanishing} naturally fits into a wider circle of ideas relating to divisors with positive normal bundles. In this direction, we can also prove the following result, which is not needed for the rest of this paper but {plays an important role in our subsequent paper \cite{Conlon4}}.

\begin{prop} Let $(X,D)$ be a complex orbifold pair as in Proposition \ref{t:vanishing} without the assumption that $X$ is K\"ahler or that $-K_X = q[D]$. Then $D$ is almost ample in the sense of Tian-Yau; compare \cite[Definition 1.1(ii)]{Tian}. Moreover, if $X$ is smooth and K\"ahler, then it is projective algebraic.
\end{prop}

\begin{proof}
As in the proof of Lemma \ref{heythere}, we begin by equipping $[D]$ with a Hermitian metric of strictly positive curvature in some tubular neighbourhood of $D$.
We can easily use this metric to show that the complex manifold $X \setminus D$ is $1$-convex so that we may take its Remmert reduction $V$ (e.g.~see
\cite[Appendix A]{Conlon}).
By construction, $V$ is biholomorphic to $X\setminus D$ in a neighbourhood of infinity. Thus, we may
compactify $V$ as a normal compact complex space $\oV$ by adding the orbifold divisor $D$.

Let $m \in \N$ be so large that $mD$ makes sense as a Cartier divisor on $\oV$. Then, from \cite[p.~347, Satz 4]{Grau:62} we see that the line bundle $[mD]$ associated with $mD$ is positive in the sense of Grauert \cite[p.~342, Definition 2]{Grau:62}, so that $\oV$ admits an embedding into $\P^{N}$ for some $N$ from \cite[p.~343, Satz 2]{Grau:62}. So far, we have essentially followed the proof of \cite[Lemma 2.1]{Epstein}.

It is clear from the proof of \cite[p.~343, Satz 2]{Grau:62} that the embedding $\oV \to \P^N$ is given by the global sections of $m'[mD]$ for some $m' \in \N$. In particular, it pulls back $\mathcal{O}_{\P^N}(1)$ to $m'[mD]$. Pulling back further under the Remmert reduction map $X \to \oV$, we obtain the line bundle $m'm[D]$ on $X$. This precisely says that $D \subset X$ is almost ample in the sense of \cite[Definition 1.1(ii)]{Tian}.

{Now $X$, as a proper modification of the projective variety $\overline{V}$, is a Moishezon space. Hence $X$ is projective algebraic if it is moreover smooth and K\"ahler \cite[Theorem 2.2.26]{ma1}.}\footnote{There certainly exist smooth pairs $(X,D)$ with $D$ positively embedded and $-K_X = q [D]$ for some integer $q > 1$ such that $X$ is Moishezon but not K\"ahler, or equivalently, not projective. See for instance \cite[pp.~1982--1985]{chnp}.}
\end{proof}

{Notice that once we know that $D$ is almost ample, we can reprove Lemma \ref{heythere} by pulling back the Fubini-Study metric on $\mathcal{O}_{\P^N}(1)$; but such a proof would rely on Grauert's deep work \cite{Grau:62}.}

\subsection{Completion of the proof of Theorem \ref{mainmain} when $p = 1$}\label{complete}
As explained at the beginning of Section \ref{hippers}, it remains to show that all K\"ahler classes on $X \setminus D$ are $\mu$-almost compactly supported in the sense of \cite[Definition 2.3]{Conlon}. We will now prove this with $\mu = -2$, relying on Proposition \ref{t:vanishing}. Indeed, we have a slightly stronger statement, which is in fact needed to treat the case $n = 2$.

\begin{prop}\label{refer}
Let $(X,D)$ be a K\"ahler orbifold pair as in Proposition \ref{t:vanishing} with $q > 1$. Then the restriction map $H^{1,1}(X,\R) \to H^2(X\setminus D, \R)$ is surjective. In particular, every K\"ahler class on $X \setminus D$ is $(-2)$-almost compactly supported in the sense of \cite[Definition 2.3]{Conlon} with $K = \emptyset$.
\end{prop}

\begin{proof}[Proof] Cohomology here means orbifold de Rham cohomology. Since $D$ is Fano by adjunction, the orbifold Calabi-Yau theorem tells us that $D$ admits K\"ahler metrics of positive Ricci curvature, so that $H^{1}(D,\R)=0$ by the usual Bochner argument and orbifold Hodge theory. It then follows from the orbifold Gysin sequence (\ref{Gysin}) that the restriction map $H^{2}(X,\R)\to H^{2}(X\setminus D,\R)$ is surjective. But $H^2(X,\R) = H^{1,1}(X,\R)$ because $h^{2,0}(X)=h^{0,2}(X)=0$ by Proposition \ref{t:vanishing}.

Thus, for any K\"ahler form $\omega$ on $X\setminus D$, we can always find a closed real $(1,1)$-form $\xi$ on $X$ such that $\omega-\xi=d\eta$ on $X \setminus D$ for some real $1$-form $\eta$ on $X\setminus D$. Using (\ref{reg_est1}), (\ref{reg_est2}), and (\ref{elephant}) from the proof of Proposition \ref{rate}, one can easily show that $\exp^{*}(\xi)=O(r^{-2})$ with $g_{0}$-derivatives.
\end{proof}

This finishes the proof of Theorem \ref{mainmain} when $p = 1$. We close this section by giving some context
for Proposition \ref{refer}. When proving a Calabi-Yau type existence theorem on an open K\"ahler manifold $M$, one goal is to produce Ricci-flat K\"ahler forms in as many classes $\mathfrak{k} \in H^2(M,\R)$ as possible. It is necessary to assume that $\mathfrak{k}$ contains \emph{some} K\"ahler form, and in the asymptotically conical setting, our ``$\mu$-almost compactly supported'' condition on $\mathfrak{k}$ is sufficient. {Goto \cite{goto} showed that the former implies the latter (with $\mu = -2$) if $M$ is a crepant resolution of a Calabi-Yau cone; Proposition \ref{refer} accomplishes the same for compactifiable manifolds $M = X \setminus D$ as in Theorem \ref{mainmain}. Moreover, in \cite{Conlon4}, we show that the form $\xi$ in the proof of Proposition \ref{refer} can be taken to be a K\"ahler form.}

\subsection{Extension from $p = 1$ to $p \in \N$}\label{s:torsion_case} Up until this point, we have been assuming for convenience that $-K_X = q [D]$ with $\alpha = q > p = 1$. This is precisely the subcase of Theorem \ref{mainmain} where $K_M$ is not only torsion but trivial. We now explain two ways of extending our results to the general case.

\subsubsection{Branched coverings}\label{ss:branched}
The first approach is more intuitive but only works if $X$ is smooth. Recall from Section \ref{s:intro:1} that, in this case, we can assume that $p,q$ are coprime.

\begin{lemma}\label{l:cyclic}
Let $X$ be a compact complex manifold, and let
$D$ be a smooth divisor in $X$ satisfying $-pK_X = q[D]$, where $p,q \in \N$ are coprime.
Then we can construct a $p$-fold cyclic covering $\tilde{X}$ of $X$, branched only over $D$, such that  $-K_{\tilde{X}}=(1+q-p)[\tilde{D}]$ with $\tilde{D}$ the preimage of $D$ in $\tilde{X}$.
\end{lemma}

This is a standard construction \cite[I.17]{Barth}. Notice that $[D]$ is divisible by $p$ in ${\rm Pic}(X)$, as required for the construction of a branched covering, precisely because $p$ and $q$ are coprime. In the setting of Theorem \ref{mainmain}, the results of Sections \ref{s:covering}--\ref{complete} then apply to the pair $(\tilde{X}, \tilde{D})$ from Lemma \ref{l:cyclic} and to the pullback $\tilde{\mathfrak{k}}$ of the given K\"ahler class $\mathfrak{k}$. The resulting Calabi-Yau metrics $\tilde{\omega}_c$ are invariant under the deck group by the final assertion of Theorem \ref{mainmain} because $\tilde{\mathfrak{k}}$ is invariant by definition.

For a classical example, let $X = \P^2$ and $D$ $=$ smooth quadric ($p = 2$, $q = 3$). Then $\tilde{X} = \P^1 \times \P^1$, $\tilde{D}$ $=$ diagonal, $\tilde{X}\setminus\tilde{D}$ is diffeomorphic to $T^*\mathbb{S}^2$, the Ricci-flat metric on $\tilde{X} \setminus \tilde{D}$ is the Eguchi-Hanson metric, and $X \setminus D = T^*\R\P^2$ represents a free isometric $\Z_2$-quotient of Eguchi-Hanson.

\subsubsection{Multivalued volume forms} The branched covering trick may not be applicable when $X$ is not smooth. However, the following more abstract argument lets us treat all cases simultaneously. Fix a defining section $s$ of $[D]$, and consider the equation $\Omega^{\otimes p} \otimes {s^{\otimes q}} = 1$ for some unknown $(n,0)$-form $\Omega$. (Here we have made use of the relation $pK_X + q[D] = 0$.) At each point of $M$, this equation has $p$ solutions, differing from each other by the $p$-th roots of unity. Thus, analytic continuation yields a \emph{multivalued} holomorphic volume form $\Omega$ on $M$, and each branch of $\Omega$ blows up to order $\alpha = q/p$ along $D$. All of the relevant work in \cite{Conlon} and in this paper generalises to this setting; in particular, the equation $(\hat\omega_c + i\partial\bar\partial u)^n = i^{n^2}\Omega \wedge \bar\Omega$ still makes sense. {(To extend Proposition \ref{t:vanishing}, notice that,  as before, it suffices to show that $-K_X$ admits a nonnegatively curved Hermitian metric with strictly positive curvature near $D$. Since $-pK_X = q[D]$, this follows from Lemma \ref{heythere} together with the fact that, for every line bundle $L$, every Hermitian metric on $pL$ has a well-defined $p$-th root on $L$.)}

Unlike in Section \ref{ss:branched}, this argument provides no particular information about $\pi_1(M)$. However, in general, the situation really can be more complicated than in Section \ref{ss:branched}. In Appendix \ref{s:cyclic_quotient}, we will see many orbifold examples where $\Gamma = \pi_1(M)$ is still finite cyclic and the universal cover of $M$ can still be written as $\tilde{X} \setminus \tilde{D}$ with $(\tilde{X}, \tilde{D})$ satisfying the hypotheses of Theorem \ref{mainmain}, but $\Gamma$ now acts nontrivially on $\tilde{D}$ and no longer only rotates the fibres of $N_{\tilde{D}}$. \newpage
\section{{From K\"ahler-Einstein to Sasaki-Einstein}}\label{s:vanCapp}

\subsubsection*{Irregular Calabi-Yau cones and crepant resolutions} The only examples of Calabi-Yau cones known before 2004 were quasiregular, i.e.~holomorphically isometric to $(K_D^{\gamma})^\times$ with a Calabi ansatz metric for some K\"ahler-Einstein Fano orbifold $D$ and $\frac{1}{\gamma} \in \N$. In 1994, Cheeger-Tian \cite{Cheeger} even conjectured that irregular Calabi-Yau cones did not exist. This conjecture was disproved in 2004 by the explicit irregular examples of Gauntlett-Martelli-Sparks-Waldram \cite{Sparks2, Sparks}. Futaki-Ono-Wang \cite{futaki} further gave an analytic construction of a large class of irregular toric examples in 2009.

From the point of view of Cheeger-Tian \cite{Cheeger}, the question remained whether irregular cones can arise as tangent or limit cones of smooth K\"ahler-Einstein manifolds. This has since been confirmed in a number of papers. In fact, the affine algebraic variety underlying an irregular Calabi-Yau cone often admits a crepant resolution, and Ricci-flat K\"ahler metrics on such crepant resolutions can be constructed using Monge-Amp{\`e}re \cite{goto, vanC2, vanC4} as well as explicit methods  \cite{lu-pope, sparks-res, oy}.

\subsubsection*{Deforming irregular Calabi-Yau cones} A fundamentally different way of desingularising algebraic varieties is by deforming rather than resolving them, thereby changing their complex rather than their symplectic structure and replacing the singularities by Lagrangian rather than holomorphic cycles; this process is especially relevant for Fano varieties \cite{br, spotti}. It seems that the only currently known irregular Calabi-Yau cone whose underlying algebraic variety can be desingularised by deformation is $K_D^\times$ with its toric cone metric from \cite{futaki}, where $D= {\rm Bl}_{p_1,p_2}\P^2$; {indeed, this is currently the only known such example in complex dimensions $4$ and less \cite{Conlonlon}.}

Van Coevering \cite[Example 6.2]{vanC} pointed out that the unique deformation of $K_D^\times$, $D =  {\rm Bl}_{p_1,p_2}\P^2$, is isomorphic to the affine manifold $X \setminus D$ with $X = {\rm Bl}_p\P^3$ and $D \in |{-\frac{1}{2}K_X}|$ realised as the strict transform of a smooth quadric passing through $p$ in $\P^3$. One then naturally wonders whether or not a Tian-Yau approach can be set to work for this example; the standard Tian-Yau setting only allows for quasiregular cones. This was in fact claimed in \cite{vanC} as an application of \cite[Theorem 1.3]{vanC}.

\subsubsection*{Convergence of complex structures} The issue with this claim is that the proof of \cite[Theorem 1.3]{vanC} crucially relies on the $\nu$-th formal neighbourhoods $D_{(\nu)}, \tilde{D}_{(\nu)}$ of $D$ in $X, N_D$ being isomorphic for all $\nu \geq 2$; this is to ensure rapid convergence of the complex structures of $X$ and $N_D$ \cite[p.~17]{vanC}. The isomorphism
$D_{(\nu)} \cong \tilde{D}_{(\nu)}$ is obtained by induction, where the inductive step follows from a certain cohomology vanishing \cite[p.~2]{vanC} that holds in the example of interest. However, contrary to what is claimed in \cite[p.~17]{vanC}, the base step $D_{(2)} \cong \tilde{D}_{(2)}$, which by \cite[Proposition 1.5]{Abate} is equivalent to the tangent sequence $0 \to TD \to TX|_D \to N_D \to 0$ being split, does not hold in the example.

In fact, if $D$ is any smooth divisor in a compact complex manifold $X$ such that $N_D$ is positive and such that $D_{(\nu)} \cong \tilde{D}_{(\nu)}$ for all $\nu \geq 2$, then, by \cite[Satz 4]{Grauert}, $X$ and $N_D$ must already be biholomorphic in some honest tubular neighbourhood of $D$. \cite[Theorem 6.6]{Rossi2} then implies that the rings of global holomorphic functions on the Remmert reduction $V$ of $X \setminus D$ \cite[Appendix A]{Conlon} and on the normal cone $(N_D^*)^\times$ are isomorphic, so that $V \cong (N_D^*)^\times$ by \cite[Chapter V, \S 7]{theory}. In other words, $X \setminus D$ is a resolution of singularities of the cone $(N_D^*)^\times$; in particular, the statement of \cite[Theorem 1.3]{vanC} was already known \cite{goto, vanC2}. However, in the example, $X \setminus D$ is affine.\footnote{Alternatively, one can use that the strict transform of a smooth quadric through $p \in \P^3$ does not appear on the known list \cite[Theorem 6.5]{Jan} of all submanifolds of ${\rm Bl}_{p} \P^3$ whose tangent sequence splits.}

\subsubsection*{Outlook} Thus, for pairs $(X,D)$ as in Theorem \ref{mainmain} where $D$ is not K\"ahler-Einstein, yet the normal $\mathbb{S}^1$-bundle of $D$ admits an irregular Sasaki-Einstein structure inducing the given CR structure, and $X \setminus D$ is not a blow-up of $(N_D^*)^\times$, the Tian-Yau problem has not been solved in the literature.

In the rest of this paper, we solve this problem for $X = {\rm Bl}_p \P^3$ and $D = {\rm Bl}_{p_1,p_2}\P^2$, currently the only pair we know of satisfying these assumptions. Our solution depends on numerical coincidences and very specific properties of these particular manifolds, so the general problem of proving a useful ``irregular'' Tian-Yau theorem remains open. \newpage
\section{An affine Calabi-Yau manifold with irregular tangent cone}\label{s:irregular}

\subsection{Overview} We consider the Fano $3$-fold $X = {\rm Bl}_p\P^3$, with $D \in |{-\frac{1}{2}K_X}|$ the strict transform of a smooth quadric passing through $p$, and we fix one of the irregular Calabi-Yau
cone structures $(g_0,\Omega_0)$ on $N_D\setminus 0$ constructed in \cite{futaki}. Our goal is to prove Theorem \ref{t:second_main}:~the existence of AC Calabi-Yau structures $(g,\Omega)$ on the affine variety
$M = X \setminus D$ that are asymptotic to $(g_0, \Omega_0)$ at infinity. In fact, almost all of the work in this section goes into proving a somewhat weaker result.

\begin{theorem}\label{irr-exi}
Let $r$ denote the radius function of the cone metric $g_0$. There exists a diffeomorphism $\Phi$ from a neighbourhood of the zero section of $N_D$ onto a neighbourhood of $D$ in $X$
such that for all $\mathfrak{k} \in H^2(M)$ and all $c > 0$, there exists a Calabi-Yau metric $\omega_c \in \mathfrak{k}$ on $M$ such that
\begin{equation}\label{e:irr-exi}
|\nabla^{k}_{g_0}(\Phi^{*}g_{c}-c g_{0})|_{g_{0}} \leq C(k) r^{-0.6752 - 0.3376 k}
\end{equation}
for all $k \in \N_0$. If $\mathfrak{k} = 0$, then we can improve $0.6752$ slightly to $0.6816$.
\end{theorem}

While the rate for $k = 0$ in (\ref{e:irr-exi}) is faster than in (\ref{e:second_main}), the improvement gained per derivative is strictly less than $1$. This is ultimately because Theorem \ref{irr-exi} is proved by solving a Monge-Amp{\`e}re equation over a rough background. However, since $g_c$ solves the Einstein equation ${\rm Ric}(g_c) = 0$, an elliptic equation modulo diffeomorphisms, Theorem \ref{t:second_main} can be deduced by gauge fixing.

\begin{proof}[Proof of Theorem \ref{t:second_main}]
The relevant ideas were introduced in \cite[Sections 2--3]{Cheeger}; see also \cite[Section 7]{Biquard}, although the problem solved there is global rather than asymptotic, hence harder.
Set $c = 1$.

The goal is to find a complete vector field $Y$ on a neighbourhood of infinity in the cone such that
$h = \exp(Y)^*\Phi^*g_1 - g_0$ satisfies the Bianchi gauge condition $\mathfrak{B}_{g_0}(h) = {\rm div}_{g_0}(h - \frac{1}{2}{\rm tr}_{g_0}(h)g_0) = 0$ and such that $h = O_{-0.3376}(r^{-0.0128})$ (cf.~Definition \ref{d:deltadecay} below); with this, (\ref{e:second_main}) follows from the elliptic system ${\rm Ric}(g_0 + h) = 0$, $\mathfrak{B}_{g_0}(h) = 0$ by bootstrapping. Using the implicit function theorem, we can reduce the construction of $Y$ to solving the linearised problem $\Delta_{g_0}Y + \mathfrak{B}_{g_0}(\Phi^*g_1 - g_0) = 0$, which can be done with $Y \in C^{2,\alpha}_{0.9872}$ \cite[Definition 2.7]{Conlon} and $\nabla^{2}_{g_0} Y = O_{-0.3376}(r^{-1.0128})$. Then $Y$ is indeed complete
and we have $h = (\Phi^*g_1 - g_0) + L_{Y}g_0$ $+$ \emph{higher order terms} $=$ $O_{-0.3376}(r^{-0.0128})$.

The uniqueness statement in Theorem \ref{t:second_main} is a direct application of \cite[Theorem 3.1]{Conlon}.

For the $(\mathbb{S}^1)^2$-invariance, we first note that $(\mathbb{S}^1)^2 \subset {\rm Aut}(X,D) \subset {\rm Aut}(D)$ also acts on the cone by holomorphic isometries. With more care in Section \ref{s:irr_estimates_better_diffeo} (using an $(\mathbb{S}^1)^2$-invariant partition of unity to glue together the identity maps in all five charts $K_i$), we can then arrange for $\Phi$ in Theorem \ref{irr-exi} to be equivariant. For each holomorphic vector field $Y$ on $M$ induced by the $(\mathbb{S}^1)^2$-action, it is then clear from (\ref{e:irr-exi}) that the exact $(1,1)$-form $\eta = L_Y\omega_c$ satisfies $|\Phi^*\eta|_{g_0} \leq C r^{-0.6752}$. Writing $\eta = i\partial\bar\partial u$ locally, we have $\Delta_{g_c}u = 0$ from linearising the Calabi-Yau condition.  Thus, $\eta$ is $g_c$-harmonic. Since $g_c$ is AC in the sense of \cite[Definition 1.11]{Conlon}, it follows that $\eta \in C^\infty_{-0.6752}(M,g_c)$ \cite[Definition 2.7]{Conlon}. Consequently, from \cite[Corollary 3.9, Theorem 3.11]{Conlon}, $\eta = 0$ as desired.
\end{proof}

In the rest of this section, we first state some useful properties of $X$ and $D$ in Section \ref{s:setup}, and then
construct suitable background metrics in Section \ref{s:irr:background} (the hardest step). Finally, we solve the Monge-Amp{\`e}re equation on these backgrounds in Section \ref{irr-MA}, thereby proving Theorem \ref{irr-exi}.

\subsection{Set-up}\label{s:setup} In this section, we collect together some of the main facts about the pair $(X,D)$ that we need for the proof of Theorem \ref{irr-exi}. We defer computational details to Appendix \ref{comp4}.

\subsubsection{Basic properties of $X$ and $D$}

The pair $(X,D)$ is a del Pezzo $3$-fold \cite[Definition 3.2.1]{AG5} with Aut$(X,D) = (\C^*)^2$. Since $D^3 = (-K_D)^2 = 7$, the linear system $|D|$ is very ample, embedding $X$ as an intersection of quadrics in $\P^8$  \cite[Proposition 3.2.4]{AG5}. On the other hand, $X$ fits into a canonical $\C^*$-family $\{X_t\}_{t \in \C}$ \cite[Proposition 5.1]{Conlon} whose central fibre $X_0$ is the $\P^1$-bundle $\mathbb{P}(K_D \oplus \C)$ with its zero section blown down.  Since this is a toric singularity, the methods of \cite{altmann} yield that $\{X_t\}_{t \in \C}$ is a versal deformation of $X_0$
and provide explicit equations for $X_t$ in $\P^8$, which can be found in \cite{dP2-physics}. As a result, $X = X_1$ is cut out by $14$ quadrics in $\P^8$, hence, in particular, is not a complete intersection. We have recorded these quadrics in (\ref{eqns-mfd}) using homogeneous coordinates $[z_1:\ldots:z_9]$ on $\mathbb{P}^8$.

With this particular choice of coordinates, the hyperplane section $X_t \cap (z_9 = 0)$ is isomorphic to $D$ and remains fixed throughout the deformation. Thus, we can derive affine equations for the cone $K_D^\times$ in $\C^8$ by deleting all monomials containing $z_9$ from (\ref{eqns-mfd}) and then substituting $z_i = Z_i z_9$ for
$i = 1, \dots, 8$. The resulting $14$ quadrics in $Z_1, \dots, Z_8$ cutting out $K_D^\times$ are given in (\ref{eqns-cone}).

Finally, we read from \cite[Proposition 5.1]{Conlon} that $b^2(X\setminus D) = 1$ (which confirms that $X \setminus D$ is
not a complete intersection) and that every element of $H^2(X \setminus D) \cong \R$ is a K\"ahler class.

\subsubsection{Affine coordinate charts}\label{4:charts:X} Define $p_9 = [0:\ldots:0:1] \in X \setminus D$. Then $X \setminus \{p_9\}$ can be written as the union of the affine opens $X_i = X \cap \{z_i \neq 0\}$ for $i \in \{1,2,3,5,8\}$. For each $i$, we have coordinate functions $(u_i, v_i, w_i)$ on $X_i$, giving rise to an isomorphism with $\C^3$. These functions are given by
$$
u_i = \frac{z_{p(i)}}{z_i},\;\, v_i = \frac{z_{q(i)}}{z_i}, \;\, w_i = \frac{z_9}{z_i}, \qquad
\begin{tabular}{c|c c c c c}
$i$ & 1 & 2 & 3 & 5 & 8 \\\hline
$p(i)$ & $2$ & $1$ & $1$ & $2$ & $6$  \\
$q(i)$ & $3$ & $5$ & $6$ & $7$ & $7$  \\
\end{tabular}.
$$
In particular, $w_i$ is a local defining function for $D \cap X_i$. One can then check that the locally defined meromorphic $(3,0)$-forms
$\sigma(i) {du_{i}\wedge dv_{i}\wedge dw_{i}}/w_i^2$, where $\sigma(i) = +1$ for $i = 1,3,8$ and $\sigma(i) = -1$ for $i = 2,5$, glue up as a meromorphic $(3,0)$-form $\Omega$ on $X \setminus \{p_9\}$ (hence on $X$) with double poles along $D$ and no other poles or zeros elsewhere.\footnote{Of course, {because $-K_X = 2[D]$}, we already knew that such a form $\Omega$ had to exist and is moreover unique up to a constant factor. The use of the argument here is to provide convenient local expressions for $\Omega$.} See Appendix \ref{charts-mfd} for more details.

As usual, the coordinates $(u_i, v_i, w_i)$ on $X$ induce coordinates $(U_{i},V_{i},W_{i})$ on the normal bundle to $D$ in $X$, where $U_i, V_i$
are simply the restrictions of $u_i, v_i$ to $D = (w_i = 0)$ in the $i$-th chart, and $W_i$ is the fibre coordinate associated with the local trivialising section $\partial_{w_i}|_{w_{i}=0}$. It is then clear that
$U_i = z_{p(i)}/z_i$, $V_i = z_{q(i)}/z_i$, and $W_i = z_{9}/z_i$ under the projective completion in $\P^8$ of the embedding (\ref{eqns-cone}) of $N_D \setminus 0$ into $\C^8$. The locally defined $3$-forms $\sigma(i){dU_{i}\wedge dV_{i}\wedge dW_{i}}/W_i^2$ again glue up as a global $3$-form $\Omega_0$ on $N_D$ with double poles along the zero section. Moreover, up to a constant factor, $\Omega_0$ is the
unique such form that is homogeneous under the standard $\C^*$-action on $N_D$.

\subsubsection{The irregular Calabi-Yau cone metric on $N_D \setminus 0$} Since $N_D = K_D^*$, \cite[Corollary 1.3]{futaki} provides us with an irregular Calabi-Yau cone metric $\omega_0$ on $N_D \setminus 0$, with infinite end at the zero section, such that $\omega_0^3 = i\Omega_0 \wedge \bar\Omega_0$. More precisely, $\omega_0 = \frac{i}{2}\partial\bar\partial r^2$, where the distance function $r$ is a smooth proper function on $N_D \setminus 0$ that diverges to $+\infty$ at the zero section. In analogy with Section \ref{s:regular:scaling}, all we need to know about $r$ here is an explicit formula for the scaling map $\nu_t = \exp((\log t) r\partial_r)$. This we can deduce from \cite{Tristan, futaki, Yau} (see in particular the last example of \cite[Section 6]{Tristan}). Indeed,
\begin{equation}\label{irr:flow}
\nu_{t}(Z_{1},\ldots,Z_{8})=(t^{a_{1}}Z_{1},\ldots,t^{a_{8}}Z_{8}),
\end{equation}
where the weights $a_i \in \R^+$ are given by
\begin{equation}\label{irr:weights}
\begin{split}
&a_1 = a_2 = \frac{9}{16}(\sqrt{33}-1) \approx 2.6688,\\
&a_3 = a_4 = a_5 = 3,\\
&a_6 = a_7 = \frac{1}{16}(105-9\sqrt{33}) \approx 3.3312,\\
&a_8 = \frac{9}{8}(9-\sqrt{33}) \approx 3.6624.
\end{split}
\end{equation}

\begin{figure}
\caption{The irregular scaling flow projected to $D$.}
\vskip5mm
\label{flow-pic}
\includegraphics{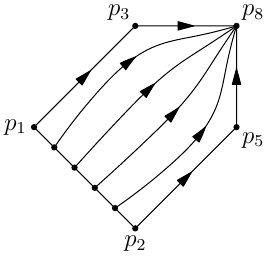}
\end{figure}

Figure \ref{flow-pic} shows the geometry of the scaling flow $\nu_t$. If the cone structure was \emph{regular} ($a_i = a$ for all $i$), then the orbits of $\nu_t$ would flow straight down to $D$ within the fibres of $N_D$.
In our \emph{irregular} setting, we have an additional ``drift'' in directions parallel to $D$, as shown in the figure; specifically, $D$ is toric, and we have drawn the projection of the orbits of $\nu_t$ to the moment polygon of $D$. The corners correspond to the fixed points of the torus action on $D$, which, in homogeneous coordinates, are simply the points $p_i = [z_1: \ldots: z_8: 0]$ with $z_j = \delta_{ij}$ for $i \in \{1,2,3,5,8\}$. Notice that $p_i = (0,0,0)$ in our chart $(U_i, V_i, W_i)$.  Notice also that almost all flowlines of $\nu_t$ converge to $p_8$ as $t \to \infty$.

\subsection{Construction of background metrics}\label{s:irr:background}

We now proceed to construct K\"ahler metrics on the complement $M = X \setminus D$ that are asymptotic to the irregular Calabi-Yau cone metric $g_0$ on $N_D \setminus 0$ in a sufficiently strong sense. The following definition will aid clarity in this construction.

\begin{definition}\label{d:deltadecay}
Let $C$ be a cone with metric $g_{0}$, radius function $r$, and Levi-Civita connection $\nabla_0$. Let $T$ be a tensor on $\{r > 1\}$. Given $\mu<0$ and $\lambda \in\mathbb{R}$, we
say that ``{$T=O_{\mu}(r^{\lambda})$ with $g_{0}$-derivatives}'', or more simply that ``$T=O_{\mu}(r^{\lambda})$'', if and only if $|\nabla^{k}_{0}T|_{g_{0}} \leq C(k) r^{\lambda+k\mu}$ on $\{r > 1\}$ for each $k\in\mathbb{N}_{0}$. In particular, $T=O_{-1}(r^{\lambda})$ with $g_{0}$-derivatives if and only if $T=O(r^{\lambda})$ with $g_{0}$-derivatives.
\end{definition}

We can now state the main result of this section.

\begin{prop}\label{p:background}
Define $\delta = 0.3376$ and $\epsilon = 0.6816$. Then there exists a compact set $K \subset M$ and a diffeomorphism $\Phi: \{r > 1\} \to M \setminus K$ such that
\begin{equation}\label{bg1e1}
\Phi^*J - J_0 = O_{-\delta}(r^{-\epsilon + \delta - 1}).
\end{equation}
Moreover, for all $\mathfrak{k} \in H^2(M)$ and $c > 0$, there exists a complete
K\"ahler metric $\hat{\omega}_c \in \mathfrak{k}$ such that
\begin{eqnarray} \label{e:background:1}
\Phi^{*}\hat{\omega}_c-c\omega_{0} = O_{-\delta}(r^{-2\delta}),\\ \label{e:background:2}
\log \Phi^*(\hat{\omega}_c^3/ic^3\Omega \wedge \bar\Omega) = O_{-\delta}(r^{-2\delta}).\,
\end{eqnarray}
\end{prop}

By adapting some arguments from \cite{Conlon}, Proposition \ref{p:background} can be derived from the following more technical statement, which will be easier to prove from first principles.

\begin{prop}\label{background2}
There exists a diffeomorphism $\Phi: \{r > 1\} \to M \setminus K$ such that
\begin{equation}\label{bg2e1}
\Phi^{*}\Omega-\Omega_{0}=O_{-\delta}(r^{-\epsilon+\delta-1}).
\end{equation}
Moreover, for all $\mathfrak{x} \in H^{1,1}(X,\R)$, there exists a closed real $(1,1)$-form $\xi \in \mathfrak{x}$ such that
\begin{equation}
\Phi^*\xi = O_{-\delta}(r^{-2\delta}).
\end{equation}
\end{prop}

Let us first see how Proposition \ref{p:background} follows from Proposition \ref{background2}.

\begin{proof}[Proof of Proposition \ref{p:background}] We need to adapt the beginning of the proof of \cite[Theorem 2.4]{Conlon} on \cite[pp.~2869--2870]{Conlon} to our setting. To begin with, (\ref{bg2e1}) implies (\ref{bg1e1}) via the same idea as in the proof of \cite[Lemma 2.14]{Conlon}. Next, we require the conclusion of \cite[Lemma 2.15]{Conlon}; the hypothesis of that lemma is satisfied here with $\lambda = -\epsilon$. Also, \cite[(2.8)]{Conlon} holds with $O(r^{-\epsilon-k\delta})$ on the right-hand side. Finally, according to the proof of Proposition \ref{refer},
$\mathfrak{k} = \mathfrak{x}|_M$ for some class $\mathfrak{x} \in H^{1,1}(X,\R)$, so that $\mathfrak{k}$ can be represented by $\xi|_M$. Since $2\delta < \epsilon$, the proposition follows.
\end{proof}

The rest of this section is concerned with the proof of Proposition \ref{background2}. This will be analogous in spirit to the proof of Proposition \ref{rate}, which dealt with the quasiregular setting. However, the details here are more complicated; in particular, it will be important to optimise the values of $\delta$ and $\epsilon$.

The first step (Section \ref{s:irregular:scaling}) is to estimate the norm of smooth functions and of their derivatives on $N_D$ with respect to the irregular cone metric. This is analogous to Section \ref{s:regular:scaling}

In the second step (Section \ref{s:irr:exp}), we use these estimates to prove bounds on $\Phi^*\Omega - \Omega_0$ and $\Phi^*\xi$ by expressing these forms in local coordinates on $N_D$ as in Sections \ref{s:reg:exp} and \ref{s:reg:asympt}. Here, $\xi$ can be any element of $\mathfrak{x} \in H^{1,1}(X,\R)$ and $\Phi$ can be any smooth map satisfying the following definition.

\begin{definition}\label{d:exp-type}
An \emph{exponential-type map} is a smooth map $\Phi$ from a neighbourhood $T$ of the zero section $D \subset N_D$ to a neighbourhood of $D \subset X$ such that $\Phi(p) = p$ for all $p \in D$ and such that $d\Phi|_p$ is complex linear with $d\Phi|_p(v) + T^{1,0}_pD = v$ for all $p \in D$ and all $v \in N_{D,p} \subset T_p^{1,0} N_D$.
\end{definition}

This definition is motivated by the constructions of Appendix \ref{app:c1}. Notice that, after shrinking $T$ if need be, $\Phi$ is necessarily a diffeomorphism onto its image.

Now, the main difference with the quasiregular case considered in Section \ref{s:asy_hol_vf} is that the resulting bounds on $\Phi^*\Omega - \Omega_0$ and $\Phi^*\xi$ are too weak to apply any of the existence results for AC Calabi-Yau metrics that we are aware of. In fact, different estimates hold in different regions, and the rates that we do obtain are sufficient \emph{only} away from some neighbourhood of the highest fixed point $p_8$.

On the other hand, since issues arise only in one particular chart, we can hope to redefine $\Phi$ and $\xi$ in a simple manner on this chart to improve the decay rate there, while preserving the exponential-type property of $\Phi$, the cohomology class of $\xi$, and the good estimates that hold everywhere else. It turns out that this is indeed possible; see Section \ref{s:irr_estimates_better_diffeo}. This will finish the proof.

\subsubsection{Estimates for smooth functions in terms of the irregular cone metric}\label{s:irregular:scaling}
The embedding (\ref{eqns-cone}) of $K_D^\times$ into $\C^8$ shows that a neighbourhood of infinity of $K_D^\times$ is covered by the five open sets
\begin{equation}\label{d:Ki}
K_{i} =\{Z \in K^{\times}_{D}: R < |Z| < 3|Z_i|\},
\end{equation}
where $i \in \{1,2,3,5,8\}$ and $R \geq 1$ is fixed as large as required; indeed, if $Z \in \C^8 \setminus \{0\}$, then obviously $|Z| < 3|Z_i|$ for some $i \in \{1,\ldots,8\}$, but if in addition $Z \in K_D^\times$ and $i \in \{4,6,7\}$, then,  from (\ref{eqns-cone}), $|Z_j| \geq |Z_i|$ for some $j \in \{3,5,8\}$. Notice also that $K_{i}$ is contained in the domain of definition of the coordinate chart $(U_i, V_i, W_i)$; in fact, we even have containment in a closed polydisk,
\begin{equation}\label{inclusion}
K_i \subset K_i' =  \{|U_{i}|\leq 3\} \cap \{|V_{i}|\leq 3\} \cap \{|W_{i}|\leq 3R^{-1}\}.
\end{equation}

\begin{prop}\label{p:irr:scaling}
Define $m =\min_i a_{i} = a_1 \approx 2.6688$ and $M =\max_i a_{i} = a_{8} \approx 3.6624 < m+1$. Fix any index $i \in \{1,2,3,5,8\}$. Then the following estimates hold on the region $K_i${\textup{:}}
\begin{align}
\label{est:Ui} dU_i &=O_{a_{i}-m-1}(r^{\max\{a_{i},a_{p(i)}\}-m-1}),\\
\label{est:Vi} dV_i &= O_{a_{i}-m-1}(r^{\max\{a_{i},a_{q(i)}\}-m-1}),\\
\label{est:Wi} W_{i} &=O_{a_{i}-m-1}(r^{-m}),\\
\label{est:ratio} \bar{W}_i/W_i &=O_{a_{i}-m-1}(1),\\
\label{est:dlogWi} dW_i/W_i &= O_{a_i-m-1}(r^{a_i-m-1}).
\end{align}
Moreover, if $A$ is any smooth function defined on the compact set $K_i'$, then, on $K_i$,
\begin{equation}\label{est:anysmooth}
A=O_{\max\{a_{i},a_{p(i)},a_{q(i)}\}-m-1}(1).
\end{equation}
\end{prop}

\begin{remark}\label{r:irr:not:invt}
Unlike in the regular setting, where we would have that $a_i = a_j$ for all $i,j$, the regions $K_i$ are not forward invariant under the scaling flow $\nu_t = \exp((\log t)r\partial_r)$ of the irregular cone metric unless $i = 8$. Hence, for all $i < 8$ and almost all $Z \in K_D ^\times \cap \mathbb{S}^{15} \subset \C^8$, the forward orbit $\{\nu_t(Z)\}_{t > 1}$ spends only a finite amount of time in $K_i$. However, this time is not bounded above independent of $Z$, so the estimates of Proposition \ref{p:irr:scaling} are not vacuous even when $i < 8$.
\end{remark}

The proof relies on two basic estimates. First, for all $i \in \{1,\ldots, 8\}$, it holds that
\begin{equation}\label{eq:irr:scaling:est}
Z_{i},\bar{Z}_{i},|Z_{i}|=O_{-1}(r^{a_{i}})
\end{equation}
on the whole region $\{r > 1\}$. This is clear from (\ref{irr:flow}) together with \cite[Lemma 1.6]{Conlon}. Second,
\begin{equation}\label{eq:irr:scaling:est:2}
c r^m \leq |Z| \leq C r^M
\end{equation}
for some $c,C \in (0,\infty)$ and again for all $\{r > 1\}$. This follows by considering the link $L = \{r = 1\}$ of the given Calabi-Yau cone structure:~a compact real hypersurface of $K_D^\times$, diffeomorphic to $K_D^\times \cap \mathbb{S}^{15}$ and contained in the annulus $\{c \leq |Z| \leq C\}$ for some constants $c, C \in (0,\infty)$. We then only need to observe that $(r^{-a_1} Z_1, \ldots, r^{-a_8}Z_8) \in L$, so that
$c^2 \leq \sum r^{-2a_i}|Z_i|^2 \leq C^2$.

\begin{proof}[Proof of Proposition \ref{p:irr:scaling}]

We begin by proving (\ref{est:Wi}). The main difficulty with this is that, unlike in (\ref{eq:irr:scaling:est}), we cannot simply use scaling to deduce that $W_i = Z_{i}^{-1} = O_{-1}(r^{-a_i})$. Indeed, if $Z \in K_i$ and $Z = \nu_r(Z^\circ)$ with $r = r(Z)$ and $Z^\circ = Z^\circ(Z) \in L$, then $W_i = r^{-a_i}W_i^\circ$, but $W_i^\circ$ will not be uniformly
bounded as a function of $Z \in K_i$; compare Remark \ref{r:irr:not:invt}. To overcome this issue, we first note
from (\ref{d:Ki}) and (\ref{eq:irr:scaling:est:2}) that $|W_i| \leq Cr^{-m}$ on $K_i$. Next, an induction on $k \geq 1$ shows that
\begin{equation}\label{induct}
\nabla_0^k W_i = \sum_{\ell = 1}^k \sum_{\substack{i_1, \ldots, i_\ell > 0\\i_1 + \cdots + i_\ell = k}} Z_i^{-1-\ell} \ast \nabla^{i_1}_0Z_i \ast \cdots \ast \nabla^{i_\ell}_0Z_i.
\end{equation}
Using (\ref{d:Ki}), (\ref{eq:irr:scaling:est:2}), and (\ref{eq:irr:scaling:est}), we can then derive the following, which proves (\ref{est:Wi}):
\begin{equation*}
|\nabla_0^k W_i| \leq C(k) \sum_{\ell = 1}^k r^{(-1-\ell)m} r^{\ell a_i - k} \leq C(k) r^{-m + k(a_i - m - 1)}.
\end{equation*}

We now prove (\ref{est:ratio}). Notice that the estimate $\bar{W}_i/W_i = \bar{W}_i Z_i = O_{a_i-m-1}(r^{a_i-m})$ follows from (\ref{est:Wi}) and (\ref{eq:irr:scaling:est}) using the product rule, but this is not optimal since we can make use of the fact that $|\bar{W}_i/W_i| = 1$. Indeed, using the more precise formula (\ref{induct}), we find that for all $k \geq 0$,
\begin{align*}
\begin{split}
|\nabla_0^k(\bar{W}_i/W_i)| &= \biggl|\bar{W}_i \ast \nabla_0^k Z_i + \sum_{j = 1}^k \biggl(\sum_{\ell = 1}^j \sum_{\substack{i_1, \ldots, i_\ell > 0\\i_1 + \cdots + i_\ell = j}} \bar{Z}_i^{-1-\ell} \ast \nabla^{i_1}_0\bar{Z}_i \ast \cdots \ast \nabla^{i_\ell}_0\bar{Z}_i\biggr) \ast \nabla_0^{k-j}Z_i\biggr|\\
&\leq C(k) \biggl(\biggl(\sum_{j=0}^{k-1} r^{-m + j(a_i-m) + a_i - k}\bigg) + r^{k(a_i - m) -k}\biggr) \leq C(k)r^{k(a_i - m - 1)}.
\end{split}
\end{align*}

(\ref{est:dlogWi}) is a direct consequence of applying (\ref{est:Wi}) and (\ref{eq:irr:scaling:est}) to the identity $dW_i/W_i = - W_i dZ_i$.

The proofs of (\ref{est:Ui}) and (\ref{est:Vi}) are also fairly straightforward and identical up to replacing $p(i)$ by $q(i)$, so we only sketch the proof of (\ref{est:Ui}). We proceed from the identity
\begin{equation}\label{induct:U}
dU_{i}=W_i(dZ_{p(i)}-U_i dZ_i).
\end{equation}
Recall that $W_i$, $dZ_i$, $dZ_{p(i)}$ and all of their derivatives are controlled from (\ref{est:Wi}) and (\ref{eq:irr:scaling:est}), and that $|U_i| \leq 3$. Thus, (\ref{induct:U}) allows for an inductive proof of (\ref{est:Ui}) by applying $\nabla_0^k$ to both sides and using the fact that $\nabla_{0}^j U_i = \nabla_{0}^{j-1}dU_i$ is controlled for all $j \in \{1,\ldots,k\}$ by the inductive hypothesis.

To prove (\ref{est:anysmooth}), we first note that $|A| \leq C$ on $K_i$ and then consider the obvious identity
\begin{equation*}
dA= \frac{\p A}{\p U_{i}}dU_{i}+\frac{\p A}{\p\bar{U}_{i}}d\bar{U}_{i}
+  \frac{\p A}{\p V_{i}}dV_{i}+\frac{\p A}{\p\bar{V}_{i}}d\bar{V}_{i}
+  \frac{\p A}{\p W_{i}}dW_{i}+\frac{\p A}{\p\bar{W}_{i}}d\bar{W}_{i}.
\end{equation*}
Since all partial derivatives of $A$ are uniformly bounded on $K_i$ as well, the required estimate of $dA$ follows by applying (\ref{est:Ui}), (\ref{est:Vi}), and (\ref{est:Wi}).
More generally, (\ref{est:anysmooth}) follows by induction.
\end{proof}

\subsubsection{Pulling back by an exponential-type map}\label{s:irr:exp}

\begin{prop}\label{p:irr:pullback}
Let $\Phi$ be an exponential-type map as in Definition \ref{d:exp-type}. Let $\Omega$ and $\Omega_0$ be the given holomorphic volume forms on $X\setminus D$ and on $N_D\setminus 0$ respectively. Let $\xi$ be a closed real $(1,1)$-form on $X$. Then we have the following estimates with respect to the irregular cone metric\textup{:}\medskip

\begin{center}
\begin{tabular}{|r|c|c|c|}
\hline
                                               & on $K_1, K_2$					& on $K_3, K_5$ 				& on $K_8$					\\\hline
$\Phi^*\Omega - \Omega_0$ & $O_{-0.6688}(r^{-2.3376})$        	& $O_{-0.3376}(r^{-1.3440})$        	& $O_{-0.0064}(r^{-0.0193})$        	\\
$\Phi^*\xi$                             & $O_{-0.6688}(r^{-1.3376})$        	& $O_{-0.3376}(r^{-0.6752})$        	& $O_{-0.0064}(r^{-0.0128})$        \\\hline

\end{tabular}\hskip1mm .
\end{center}\medskip

\noindent Recall here that $K_i$ denotes the open set defined in \eqref{d:Ki}.
\end{prop}

\begin{proof} We fix $i \in \{1,2,3,5,8\}$, define $\delta_i > 0$ by $\max\{a_i, a_{p(i)}, a_{q(i)}\} - m - 1 = -\delta_i$, and work on $K_i$. For better readability, we will drop almost every subscript $i$ in what follows.

We begin by estimating $\Phi^*\xi$. For this, we note that $\Phi^*\xi$ can be written as a linear combination of wedge products of two of the basic forms $dU$, $d\bar{U}$, $dV$, $d\bar{V}$, $dW$, $d\bar{W}$, the coefficients being smooth functions defined on a neighbourhood of the zero section. It is then clear from (\ref{est:Ui}), (\ref{est:Vi}), (\ref{est:Wi}), and (\ref{est:anysmooth}) that
$\Phi^*\xi = O_{-\delta}(r^{-2\delta})$,
and the claimed estimates follow from this using (\ref{irr:weights}).

The estimate of $\Phi^*\Omega - \Omega_0$ is more involved. Recall the results of Section \ref{4:charts:X}.
Then recall from Appendix \ref{app:c1} that there exist smooth local complex coordinates $(u', v', w')$ on $X$, with $w' = w$ and the difference $(u', v') - (u, v)$ vanishing at the divisor, such that, on one hand,
\begin{equation*}
du \wedge dv \wedge dw = du' \wedge dv' \wedge dw' + w' \Upsilon \wedge dw',
\end{equation*}
for some smooth complex-valued $2$-form $\Upsilon$, so that, in particular,
\begin{equation}\label{irr:IandII}
\Phi^*\Omega=\frac{\Phi^{*}(du')\wedge \Phi^*(dv') \wedge \Phi^*(dw')}{(\Phi^{*}w')^{2}} + \frac{\Phi^{*}\Upsilon \wedge \Phi^{*}(dw')}{\Phi^{*}w'} = ({\rm I}) + ({\rm II}).
\end{equation}
On the other hand, there exist smooth functions $A_j, B_j, C_j$ ($j = 1,2,3$) such that
\begin{equation*}
\begin{split}
\Phi^{*}u'-U &=A_{1}W^2 +A_{2}W\bar{W} +A_{3}\bar{W}^{2},\\
\Phi^{*}v'-V &=B_{1}W^2 +B_{2}W\bar{W} +B_{3}\bar{W}^{2},\\
\Phi^{*}w'-W &=C_{1}W^2 +C_{2}W\bar{W} +C_{3}\bar{W}^{2},
\end{split}
\end{equation*}
so that, using (\ref{est:Wi}), (\ref{est:ratio}), (\ref{est:dlogWi}), and (\ref{est:anysmooth}),
\begin{equation}\label{irr:exp:quadratic2}
\begin{split}
\Phi^*(du') &= dU + O_{-\delta}(r^{-2m-\delta}),\\
\Phi^*(dv') &= dV + O_{-\delta}(r^{-2m-\delta}),\\
\Phi^*w' &= W(1 + O_{-\delta}(r^{-m})),\\
\Phi^*(dw') &= dW + W^2 O_{-\delta}(r^{-\delta}).
\end{split}
\end{equation}

We now combine (\ref{irr:IandII}) and (\ref{irr:exp:quadratic2}) to prove the desired estimate for $\Phi^*\Omega - \Omega_0$. Using the bounds on $dU, dV$
from (\ref{est:Ui}), (\ref{est:Vi}), together with the fact that $dW/W^2 = -dZ = O_{-1}(r^{a_i-1})$ from (\ref{eq:irr:scaling:est}),
a lengthy but completely straightforward computation yields that
$${\rm (I)} = \Omega_0 + O_{-\delta}(r^{\max\{a_{i},a_{p(i)}\}+\max\{a_{i},a_{q(i)}\}- 2m - 2 - \delta}).$$
Indeed, the relevant error term is $dU \wedge dV \wedge [\Phi^*(dw')- dW]/W^2$, and one checks, using (\ref{irr:weights}), that all other contributions are of lower order. Regarding (II), we first note that the argument used to estimate $\Phi^*\xi$ above applies verbatim to the smooth $2$-form $\Phi^*\Upsilon$, so that $\Phi^*\Upsilon = O_{-\delta}(r^{-2\delta})$. Given this, the bound on $dW/W$ from (\ref{est:dlogWi}), and (\ref{irr:weights}), one then quickly finds that
$${\rm (II)} = O_{-\delta}(r^{a_i - m - 1 - 2\delta}).$$
Using (\ref{irr:weights}), one checks that $\max\{a_i, a_{p(i)}\} + \max\{a_i,a_{q(i)}\} = a_i + \max\{a_i, a_{p(i)}, a_{q(i)}\}$, which implies that our bounds for (I) and (II) are of the same order. The claimed values follow from this.
\end{proof}

\begin{remark}
It is instructive to compare the estimates for $\Phi^*\Omega - \Omega_0$ appearing in this proof with the corresponding ones in the quasiregular case (Section \ref{s:reg:asympt}, $\alpha = 2$). In both cases, our estimates for (I) and (II) are of the same order. In the quasiregular case, the overall bound $O_{-1}(r^{-n})$ is the same as the one for a smooth $n$-form defined on a neighbourhood of the zero section; here, we can do slightly better than this ($O_{-\delta}(r^{a_i - m - 1 - 2\delta})$ vs.~$O_{-\delta}(r^{-3\delta})$). This will be crucial later on.

\end{remark}

\subsubsection{A cut-off trick}\label{s:irr_estimates_better_diffeo} Let us rewrite the rates of Proposition \ref{p:irr:pullback} as $\Phi^*\Omega - \Omega_0 = O_{-\delta_i}(r^{-\epsilon_i + \delta_i - 1})$
and $\Phi^*\xi = O_{-\delta_i}(r^{-2\delta_i})$ on $K_i$. In order for our construction of AC Calabi-Yau metrics to go through right away, we would need that the following two inequalities hold for all $i$:
\begin{itemize}
\item[(1)] $\epsilon_i > 0$,
\item[(2)] $\min\{2\delta_i, \epsilon_i\} + \delta_i > 1$.
\end{itemize}
Here (1) is needed to be able to construct any asymptotically conical background K\"ahler metrics on $X\setminus D$ at all; cf.~our derivation of Proposition \ref{p:background} from Proposition \ref{background2}, based on \cite[Lemma 2.15]{Conlon}. Given this, (2) is then needed for the PDE analysis in Section \ref{irr-MA}.

(1) and (2) are indeed satisfied for $i = 1,2,3,5$ (albeit with a very small margin in (2) for $i = 3,5$), but they both fail for $i = 8$. On the other hand, there is considerable freedom in choosing $\Phi$ (under the exponential-type condition) and $\xi$ (preserving the class $\mathfrak{x} \in H^{1,1}(X,\R)$ represented by $\xi$). The aim of this short section is to exploit this freedom to arrange that $\Phi^*\Omega - \Omega_0 = 0$ and $\Phi^*\xi = 0$ on $K_8$. It is then clear that Proposition \ref{background2} holds with $\delta = \delta_5$ and $\epsilon = \epsilon_5$, as claimed.

Thus, let us fix an exponential-type map $\Phi$, defined on some closed tubular neighbourhood $T$ of the zero section of $N_D$, and a closed real $(1,1)$-form $\xi \in \mathfrak{x}$. Let ${\chi}_0: \C^3 \to \R$ be a smooth function with ${\chi}_0 \equiv 1$ on $\{|x| \leq 3, |y| \leq 3, |z| \leq 3\}$ and ${\rm supp}({\chi}_0) \subset \{|x| \leq 4, |y| \leq 4, |z| \leq 4\}$ and define
\begin{align*}
\chi = \chi_0 \circ (U_8, V_8, 0) \;\,{\rm on}\;\, P = \{|U_8| \leq 4, |V_8| \leq 4\} \subset N_D,\\
\hat{\chi} = \chi_0 \circ (u_8, v_8, w_8)  \;\,{\rm on}\;\, X_8 = X \cap \{z_8 \neq 0\},
\end{align*}
extending $\hat{\chi}$ by zero to the whole of $X$. Shrinking $T$ if need be, we can assume that $\Phi(T \cap P) \subset X_8$. We then define $\Phi': T \to X$ by setting $\Phi' = \Phi$ on $T\setminus P$ and $\Phi' = \chi {\rm id} + (1-\chi)\Phi$ on $T \cap P$, where the latter formula is to be evaluated working in the coordinate charts $(U_8, V_8, W_8)$ and $(u_8, v_8, w_8)$. Moreover, since $X_8 \cong \C^3$, we can write $\xi|_{X_8}= i\partial\bar\partial u$ for some smooth potential $u: X_8 \to \R$,
and we define $\xi' = \xi - i\partial\bar\partial(\hat{\chi} u)$. Then clearly $(\Phi')^*\Omega - \Omega_0 = 0$ and $(\Phi')^*\xi'= 0$ on $K_8$, and all that remains to be checked, after shrinking $T$ further if necessary, is that $\Phi'$ satisfies Definition \ref{d:exp-type}.

The properties of Definition \ref{d:exp-type} only need to be verified for $p \in D \cap P$. Now simply observe that $d\Phi'|_p = \chi {\rm id} + (1-\chi)d\Phi|_p$ for all such $p$, again working in the charts $(U_8, V_8, W_8)$ and $(u_8, v_8, w_8)$. Thus, the required properties hold for $\Phi'$ because they hold for ${\id}$ (since $(U_8, V_8, W_8)$ represents the vector $W_8(\partial_{w_8} + T_p^{1,0}D) \in N_{D,p}$ at $p = (u_8, v_8, 0) = (U_8, V_8, 0)$) as well as for $\Phi$.

This completes the proof of Proposition \ref{background2}.

\subsection{Solving the Monge-Amp{\`e}re equation}\label{irr-MA}

Proposition \ref{p:background} yields background K\"ahler metrics on $M = X \setminus D$ that are asymptotically conical and whose Ricci potential decays to zero at infinity.
In order to prove Theorem \ref{irr-exi}, it remains to solve the complex Monge-Amp{\`e}re equation
\begin{align}\label{e:final_MA}
\begin{split}
(\hat{\omega}_c + i\partial\bar\partial u)^3 = ic^3\Omega \wedge \bar\Omega = e^{{f}}\hat{\omega}_c^3,\\
\Phi^{*}\hat{\omega}_c-c\omega_{0} = O_{-\delta}(r^{-2\delta}), \;\, \Phi^*{f} = O_{-\delta}(r^{-2\delta}), \;\,\delta = 0.3376,
\end{split}
\end{align}
where here, and in the rest of this section, the $O$ includes all derivatives as in Definition \ref{d:deltadecay}.

The existence theory of \cite{Conlon} does not apply to this equation because neither $\Phi^*\hat\omega_c - c\omega_0$ nor $\Phi^*{f}$ are
$O_{-1}(r^\lambda)$ for any $\lambda < 0$. However, by combining some ideas from \cite{Conlon} and the earlier literature, we will nevertheless be able to construct a solution $u$ with $\Phi^*u = O_{-\delta}(r^{\lambda})$ for every $\lambda > -4$. The main step is the following iteration lemma, which is similar to \cite[Lemma 2.12]{Conlon}.

\begin{prop}\label{p:irr_iterate}
Fix real numbers $\kappa, \lambda, \mu$ with $-6 < \lambda \leq \kappa < 0$, $-1 \leq \mu < 0$, and $\kappa + \mu + 1 < 0$. Suppose that $\Phi^*J - J_0 = O_\mu(r^{\kappa - \mu - 1})$. If there exists a
K\"ahler metric $\omega \in \mathfrak{k}$ such that
\begin{align}
\Phi^{*}{\omega}-c\omega_{0}&=O_{\mu}(r^{\kappa}),\\
\log \Phi^{*}(\omega^3/ic^3\Omega \wedge \bar\Omega) &=O_{\mu}(r^\lambda),
\end{align}
then there exists another K\"ahler metric $\omega_\sharp \in \mathfrak{k}$ such that
\begin{align}
\Phi^{*}{\omega}_{\sharp}-c\omega_{0}&=O_{\mu}(r^{\kappa}),\\
\log \Phi^{*}(\omega_\sharp ^3/ic^3\Omega \wedge \bar\Omega) &=O_{\mu}(r^{\lambda + \kappa + \mu + 1}).
\end{align}
\end{prop}

\begin{remark}
In our application, $\kappa + \mu + 1 < 0$ holds because $3\delta = 1.0128 > 1$.
\end{remark}

\begin{proof}[Proof of Proposition \ref{p:irr_iterate}]
Let $g_*$ be a  metric on $M$ extending $\Phi_*g_0$, and let $\rho \geq 1$ be a function on $M$ extending $\Phi_*r$. We write $T = O_\mu(\rho^\lambda)$ if and only if $\Phi^*T = O_\mu(r^\lambda)$. Set ${f} = \log(i c^3 \Omega \wedge \bar\Omega / \omega^3)$.

Integrating $\nabla_{g_*}f$ along $g_*$-geodesics yields that $f \in C^{0,\alpha}_{\lambda}(M,g_{*})$ for all $\alpha < |\mu|$, as defined in \cite[Section 2.2]{Conlon}.
It then follows from a standard result \cite[Theorem 2.11]{Conlon} that we can solve $\Delta_{g_{*}}u =2f$ with $u \in C^{2,\alpha}_{\lambda + 2}(M, g_*)$.
Regarding higher derivatives of $u$, it holds for all $k \in \N$ that
$$|\nabla^{\ell} f| \leq C(\ell)\rho^{\lambda + \ell\mu}\; (\forall\ell \leq k+1) \Longrightarrow f\in C^{k,\alpha}_{\lambda+ k \mu + k}(M) \Longrightarrow |\nabla^{k + 2}u| \leq C(k) \rho^{\lambda + k\mu},$$
where all metric operations are the ones associated with $g_*$. Indeed, after
integrating $\nabla^{k+1}_{g_*}f$ along geodesics to bound the $(k + \alpha)$-seminorm of $f$, the first implication becomes trivial. For the second one, we apply \cite[Theorem 2.11]{Conlon} to solve $\Delta_{g_*} \bar u = 2f$ with
$\bar u \in C^{k+2,\alpha}_{\lambda + k\mu + k + 2}(M,g_*)$. Then $\bar{u}$ satisfies the required derivative bound,  and $u-\bar{u}$ is harmonic, hence satisfies the same bound.

Using the formula $2i\partial\bar\partial = dd^c$ and our assumed estimate on $\Phi^*J - J_0$, we can find, as in the proof of \cite[Lemma 2.12]{Conlon}, a function $\chi \in C^\infty(M)$ with $1-\chi \in C^\infty_0(M)$ such that ${\omega}_{\sharp}={\omega}+i\partial\bar{\partial}(\chi u)$ is positive definite on $M$. Then, using the bound on $\Phi^*J - J_0$ once again,
\begin{equation*}
\Phi^*\omega_\sharp - c\omega_0=(\Phi^{*}{\omega}-c\omega_{0})+\Phi^*(i\partial\bar{\partial}(\chi u)) =O_{\mu}(r^{\kappa}).
\end{equation*}
It remains to estimate the new Ricci potential $f_\sharp = \log(i c^3 \Omega \wedge \bar\Omega / \omega_\sharp^3)$. Clearly
\begin{equation*}
f_{\sharp} =f -\log\biggl(1+f + \frac{1}{2}(\Delta_{g}(\chi u)-\Delta_{g_{*}}u)+\sum^{n}_{k=2}{n\choose k}\frac{{\omega}^{n-k}\wedge(i\partial\bar{\partial}(\chi u))^{k}}
{{\omega}^n}\biggr).
\end{equation*}
Since the higher order terms are $O_\mu(\rho^{2\lambda})$, hence $O_\mu(\rho^{\lambda + \kappa + \mu + 1})$, it suffices to note that
\begin{equation*}
\Delta_{g}u-\Delta_{g_{\ast}}u =({g}-g_{*}) \ast \nabla^{2} u +
g \ast \nabla g \ast \nabla u = O_\mu(\rho^{\lambda+\kappa+\mu+1}),
\end{equation*}
the covariant derivatives being taken with respect to $g_*$. Notice that we need to use our assumption on $\Phi^*J - J_0$ here as well, in order to control $\Phi^*g - g_0$ in terms of $\Phi^*\omega - \omega_0$.
\end{proof}

Proposition \ref{p:irr_iterate} allows us to assume that the Ricci potential of our background metric is $o(r^{-2})$.
The most natural approach
to solving the Monge-Amp{\`e}re equation might then be to define weighted H\"older spaces consisting of $O_\mu(r^\lambda)$ type functions
and extend the standard theory \cite[Section 2.3]{Conlon} to this setting. In fact, we have essentially just seen how $\Delta^{-1}$ acts on such spaces. However, for the sake of brevity, we will instead assemble a solution using some arguments from the literature.

\begin{proof}[Proof of Theorem \ref{irr-exi}]
Our aim is to solve (\ref{e:final_MA}). We note that Proposition \ref{p:irr_iterate} applies here with $\kappa = \max\{-2\delta,-\epsilon\} = -2\delta = -0.6752$ (see Proposition \ref{p:background}) and $\mu = -\delta = -0.3376$. Setting $\lambda = \kappa$ and $\omega = \hat\omega_c$ initially, we can achieve that $\lambda < -6$ after finitely many iterations. More precisely, this means that we can assume without loss of generality that $\hat\omega_c$ in (\ref{e:final_MA})
satisfies
\begin{align*}
\begin{split}
\Phi^{*}\hat{\omega}_c-c\omega_{0} &= O_{-0.3376}(r^{-0.6752}),\\
\Phi^*J - J_0 &= O_{-0.3376}(r^{-1.3440}),\\
\log \Phi^*(\hat{\omega}_c^3/ic^3\Omega \wedge \bar\Omega) &= O_{-0.3376}(r^{-6.0128}).
\end{split}
\end{align*}
Notice in particular that $(M, \hat g_c)$ has Euclidean volume growth and bounded geometry, its curvature tensor and all of its covariant derivatives tend to zero at infinity, and our given irregular Calabi-Yau cone is the (only) tangent cone of $(M,\hat g_c)$ at infinity in the $C^\infty$ Cheeger-Gromov sense.

Using one of the Tian-Yau theorems \cite[Proposition 4.1]{Tian} (see \cite[Proposition 3.1]{Hein} for the precise statement that we need and for an exposition of the proof), we can therefore assert that (\ref{e:final_MA}) has a solution $u \in C^\infty(M)$ such that $\sup_M |\nabla^k u| < \infty$ with respect to $\hat g_c$ for all $k \in \N_0$.

It remains to show that $u$ decays at infinity. As explained in \cite[p.~26]{Hein}, to obtain $C^0$ decay in a setting such as ours, we can either employ barriers \cite{Santoro, Tian} or Moser iteration with weights \cite{Joyce}. Both methods give that $|u| \leq C(\lambda)\rho^{\lambda}$ for all $\lambda > -4$. Applying Schauder estimates on unit-size geodesic balls in $(M, \hat g_c)$ then shows that $|\nabla^k u| \leq C(\lambda, k)\rho^{\lambda}$ for all $\lambda > -4$ and $k \in \N_0$. We can do even better by working on balls up to size $\rho^{0.3376}$:~\cite[Lemma 3.7, Proposition 3.8(ii), ``and in fact slightly more'' on p.~25]{Hein} with $\lambda_{\textup{\cite{Hein}}}= 0.3376$ yields that $u = O_{-0.3376}(\rho^{\lambda})$ for all $\lambda > -4$.
\end{proof} \newpage
\appendix
\section{Asymptotically conical Ricci-flat K\"ahler surfaces}\label{s:kronheimer}

AC Ricci-flat K\"ahler manifolds of complex dimension $n = 2$ are completely classified \cite{Kronheimer2, suvaina, wright}; they are precisely the Kronheimer ALE spaces \cite{Kronheimer} and certain quotients of Kronheimer spaces of type $A$ by free holomorphic isometric actions of finite cyclic groups. The paper \cite{rasdeaconu} essentially shows that all of these spaces can be constructed by the Tian-Yau method, although Theorem \ref{mainmain} is needed to get a definitive result. In this appendix, we summarise and clarify {the results of \cite{rasdeaconu}}.

\subsection{Kronheimer spaces}\label{s:B1} This is the honest Calabi-Yau case, where $K_M$ is trivial and (as it turns out, equivalently) $\pi_1(M) = 0$. The asymptotic cone is $\C^2/\Gamma$ with $\Gamma$ a finite subgroup of SU$(2)$ acting freely on $\mathbb{S}^3 \subset \C^2$, and all such groups appear on the list. Let $\rho: {\rm SU}(2) \to {\rm SO}(3)$ denote the usual double covering. Then the orbifold divisor $D$ at infinity in any compactification $X$ of $M$ satisfying the hypotheses of Theorem \ref{mainmain} must be the spherical orbifold $\mathbb{S}^2/\rho(\Gamma)$ (which is, of course, isomorphic to $\P^1$ as a variety, but not as an orbifold unless $\rho(\Gamma) = \{1\}$). Notice that $\rho|_\Gamma$ is an isomorphism onto its image if and only if $-{\rm id}_{\C^2} \not\in \Gamma$, and this holds if and only if $\Gamma$ is of type $A_{k-1}$ with $k$ odd.

We list all possibilities for $\Gamma$ below. As a complex manifold, $M$ is obtained by realising $\C^2/\Gamma$ as an affine surface $f_\Gamma(x,y,z) = 0$ in $\C^3$, deforming this surface, and resolving any remaining singularities. In \cite{rasdeaconu}, orbifold pairs $(X,D)$ with $M = X \setminus D$ are obtained by embedding $\C^3$ into $\P^3(a,b,c,1)$ and then closing $M$. However, for each $\Gamma$, there is only one triple $(a,b,c)$ (corresponding to the $\C^*$-action on $\C^2/\Gamma$ induced by the scaling field $r\partial_r$ on $\C^2$) such that $(X,D)$ satisfies the hypotheses of Theorem \ref{mainmain}
(the most restrictive one being that $D$ admits a constant curvature orbifold metric).

The table also lists the singularities of $X$, where $\frac{1}{i}(\ell, m)$ means $\C^2/G$ with $G = \langle {\rm diag}(\eta^\ell, \eta^m) \rangle$ for some primitive $i$-th root of unity $\eta$, and the integer $q > 1$ is determined by $-K_X = q[D]$. Then the singularities of $D$ are given by $\frac{1}{i}(\ell)$, and the cone $\C^2/\Gamma$ can also be realised as $(\frac{1}{q-1}K_D)^\times$.

Notice that $X$ is smooth only for the types $A_0$ ($M = \C^2 = \P^2 \setminus {\rm line}$) and $A_1$ ($M = {\rm Eguchi}$-${\rm Hanson}$ $=$
$(\P^1 \times \P^1) \setminus {\rm diagonal}$), and that $X$ will have singularities that are not ordinary double points unless $\Gamma = A_0, A_1, A_2, A_3, D_4$. We also mention that the $A_{k-1}$ case is discussed in \cite[p.~30]{vanC}.\medskip

\begin{center}
\begin{tabular}{|c|c|c|c|c|c|}\hline
type 	of $\Gamma$		& $|\Gamma|$ & $f_\Gamma(x,y,z)$ 		& $(a,b,c)$					& singularities of $X$ 							& $q$		\\\hline\hline
$A_{k-1}$ ($k$ odd) 	& $k$ 		& $xy + z^k$ 				& $(k,k,2)$ 					& $\frac{1}{k}(2,1), \frac{1}{k}(2,1)$ 				& $3$ 		\\\hline
$A_{k-1}$ ($k$ even) 	& $k$ 		& $xy + z^k$ 				& $(\frac{k}{2},\frac{k}{2}, 1)$		& $\frac{2}{k}(1,1), \frac{2}{k}(1,1)$ 				& $2$ 		\\\hline
$D_{k+2}$ ($k \geq 2$)	& $4k$		& $x^2y + y^{k+1} + z^2$ 	& $(k,2,k+1)$ 					& $\frac{1}{2}(1,1), \frac{1}{2}(1,1), \frac{1}{k}(1,1)$	& $2$ 		\\\hline
$E_6$ 				& $24$ 		& $x^4 + y^3 + z^2$ 		& $(3,4,6)$ 					& $\frac{1}{2}(1,1), \frac{1}{3}(1,1), \frac{1}{3}(1,1)$ 	& $2$ 		\\\hline
$E_7$ 				& $48$ 		& $x^3y + y^3 + z^2$		& $(4,6,9)$ 					& $\frac{1}{2}(1,1), \frac{1}{3}(1,1), \frac{1}{4}(1,1)$ 	& $2$ 		\\\hline
$E_8$				& $120$ 		& $x^5 + y^3 + z^2$ 		& $(6,10,15)$ 					& $\frac{1}{2}(1,1), \frac{1}{3}(1,1), \frac{1}{5}(1,1)$ 	& $2$ 		\\\hline
\end{tabular}
\end{center}\smallskip

\begin{remark}
The simplest Kronheimer spaces are those for which, as a complex manifold, $M$ is a crepant resolution of the cone $\C^2/\Gamma = (\frac{1}{q-1}K_D)^\times$. If $q = 2$, then this cone also admits a canonical \emph{partial} crepant resolution, given by the total space of $K_D$. This space carries a complete Ricci-flat K\"ahler metric provided by the Calabi ansatz, but it still has orbifold singularities:~one cyclic SU$(2)$ singularity $\frac{1}{i}(\ell, -\ell)$ {at the zero section of $K_{D}$} for each singularity $\frac{1}{i}(\ell,m)$ of $X$ at infinity.
\end{remark}

\subsection{Free quotients of Kronheimer spaces}\label{s:cyclic_quotient} By \cite{suvaina, wright}, a nontrivial free isometric group action on a Kronheimer space is necessarily holomorphic, and all such actions are given by a cyclic group $\Z_n$ ($n > 1$) acting on a Kronheimer space $M$ of type $A_{nd-1}$ ($d \in \N$) in the following manner.

As a complex manifold, $M$ is a crepant resolution of an affine surface $M_0: xy + F(z^n) = 0$ (which is smooth for generic choices of the polynomial $F$) with $F(Z) = Z^d + \textup{\textit{lower order}}$. The $\Z_n$-action on $M$ is induced from the $\Z_n$-action on $M_0$ given by $(\zeta, (x,y,z)) \mapsto (\zeta x, \zeta^{-1}y, \zeta^m z)$ for all
$n$-th roots of unity $\zeta$ and some $m \in \N$ coprime to $n$. This map acts on $K_{M_0}$ as multiplication by $\zeta^m$.

In \cite{rasdeaconu}, orbifold compactifications $(X,D)$ of $M/\Z_n$ are obtained by choosing weights $(a,b,c)$ such that $M_0/\Z_n \subset \C^3/\frac{1}{n}(a,b,c)$, embedding $\C^3/\frac{1}{n}(a,b,c)$ into $\P^3(a,b,c,n)$, and closing $M/\Z_n$. However, $D$ does not admit a constant curvature orbifold metric unless $a = b$, or equivalently, $n = 2$, $d$ odd (with $a = b = d$, $c = 1$), because this is the only case where the $\Z_n$-action on the asymptotic cone $\C^2/\Z_{nd}$ of $M$ is induced by the standard
$\C^*$-action on $\C^2$. Thus, a Tian-Yau type theorem can be applied to the orbifold pair $(X,D)$ if and only if $n = 2$ and $d$ is odd (including the case $n = 2$, $d = 1$, which leads to the $\Z_2$-quotient of Eguchi-Hanson mentioned in Section \ref{ss:branched}).

However, in order to obtain pairs $(X,D)$ such that Theorem \ref{mainmain} applies in all cases, we only need to compactify the $A_{nd-1}$-space $M$ by a pair $(\tilde{X}, \tilde{D})$ as in Section \ref{s:B1} and then observe that the $\Z_n$- action on $M$ extends to $\tilde{X}$, preserving $\tilde{D}$, in such a way that the quotient $(X,D)$ of $(\tilde{X}, \tilde{D})$ by
this extended action again satisfies all of the hypotheses of Theorem \ref{mainmain}. Let us now describe this process more explicitly. We write $k = nd$ and let $\zeta$ denote a fixed primitive $n$-th root of unity.\medskip

\noindent {\bf Case 1:~$k$ is odd.} Both of the singularities of $\tilde{X}$ take the form $\C^2_{zw}/G$, where $G = \langle {\rm diag}(\eta^2, \eta) \rangle$ for some primitive $k$-th root of unity $\eta$, and where we can assume that the holomorphic volume form of $M$ is given by $w^{-3}dz \wedge dw$. This is proved in \cite{rasdeaconu}, and by following the constructions of \cite{rasdeaconu}, we can describe how $\Z_n$ acts in this picture. In fact, the $G$-invariant rational functions on $\C^2$ are generated by $(z^k, w^k, \frac{z}{w^2})$, and $\zeta$ acts on this triple by ${\rm diag}(\zeta^{-2}, \zeta^{-1}, \zeta^m)$. Thus, $X$ again has two singularities, both of the form $\frac{1}{nk}(2 - mk,1)$. The natural way to see this is to consider the $\Z_{nk}$-action on $\C^2_{zw}$ given by $\langle {\rm diag}(\xi^{2-mk}, \xi)\rangle$, where $\xi$ is a primitive $nk$-th root of unity such that $\xi^k = \zeta^{-1}$.

We also read from this picture that the holomorphic volume form of $M$ pushes down to an $n$-fold multivalued holomorphic volume form on $M/\Z_n$, each of whose branches blows up to order $3$ along the compactifying divisor $D$. Thus, $-pK_X = q[D]$ with $p = n$ and $q = 3n$ (so that $\alpha = \frac{q}{p} = 3)$. However, we are not able to cancel any divisors of $n$ in this relation because the $n$-fold multivalued volume form on $M/\Z_n$ is not $n'$-fold multivalued for any $n' < n$.\medskip

\noindent {\bf Case 2:~$k$ is even.} Here $G = \langle {\rm diag}(\eta,\eta)\rangle$ for some primitive $\frac{k}{2}$-th root of unity $\eta$, the volume form of $M$ can be written as $w^{-2}dz \wedge dw$, and the $G$-invariant rational functions on $\C^2$ are generated by the triple $(z^{\frac{k}{2}}, w^{\frac{k}{2}}, \frac{z}{w})$. The form of the action of $\zeta$ on this triple now depends on the parity of $d$:~it is given by ${\rm diag}(\zeta^{-1}, \zeta^{-1}, \zeta^m)$ if $d$ is even, and by ${\rm diag}(-\zeta^{-1}, \zeta^{-1}, \zeta^m)$ if $d$ is odd. In both cases, by using the same method as in Case 1 above, we deduce that the singularities of $X$ can be written as $\frac{2}{nk}(1-m\frac{k}{2},1)$. However, there are now two markedly distinct possibilities.

$\bullet$ \emph{$\Z_n$ acts nontrivially on the divisor $\tilde{D}$}: This is always the case unless $n = 2$ and $d$ is odd. Then $X$ has exactly two singular points, and $p = n$, $q = 2n$, $\alpha = 2$, in complete analogy with Case 1.

$\bullet$ \emph{$\Z_n$ acts trivially on the divisor $\tilde{D}$}: This is the case if and only if $n = 2$ and $d$ is odd. If we take the quotient $\tilde{X}/\Z_n$ in the orbifold category, then $X$ will have a divisorial singularity with angle $\frac{2\pi}{n}$ along $D$. However, we need to discard such singularities when applying Theorem \ref{mainmain} and view $X$ only as a variety. As such, $X$ has at most two singular points, both of them given by $\frac{1}{2d}(1-d,2)$. Moreover, $w^2$ is a local defining function for $D$ in $X$, so that $p = 2$, $q = 3$, $\alpha = \frac{3}{2}$.\medskip

In the examples with $n = 2$ and $d$ odd, $\tilde{X}$ can perhaps be viewed as an ``orbifold branched cover'' of $X$, ramified along the suborbifold $\tilde{D}$. These are precisely the examples where the construction of \cite{rasdeaconu} produces compactifications to which Theorem \ref{mainmain} can be applied. We also observe that {none} of the orbifolds $X$ in this section are actually smooth, \emph{except} for the $n = 2$, $d = 1$ example.

\begin{remark}\label{r:obata}
Returning to Corollary \ref{c:rate}, it seems worth pointing out that we have now seen very explicitly that $\alpha \leq 3$ in complex dimension $2$, even in the orbifold case, and that equality can occur for flat cones $\C^2/\Gamma$ with $\Gamma \neq \{1\}$. In fact, $\alpha = 3$ if and only if, as in Case 1 above, $\Gamma = \Z_{n^2d} \subset {\rm U}(2)$ with both $n$ and $d$ odd (including the case $n = 1$ with $\Gamma \subset {\rm SU}(2)$ from Appendix \ref{s:B1}).
\end{remark} \newpage
\section{Technical constructions}

\subsection{Coordinates near a divisor}\label{app:c1}
Let $X$ be a complex manifold and let $D \subset X$ be a smooth divisor. Everything we say here applies with obvious changes to the orbifolds of Section \ref{s:optimal_TY}.

Let $(z_1, \dots ,z_n)$ be local holomorphic coordinates on $X$, with $z_n$ a local defining function for $D$. Then we have canonically associated local holomorphic coordinates $(w_1, \dots, w_{n})$ on the total space of the normal bundle $N_D = T^{1,0}X|_D/T^{1,0}D$, corresponding to the coset $w_n(\partial_{z_n} + T^{1,0}D)$ of normal vectors based at $(w_1, \dots, w_{n-1},0)$. In this section, we wish to discuss a useful degree of freedom to modify the coordinates $z_i$ without changing the associated coordinates $w_i$.

For this, set $z'_i = z_i - A_i(z_1,\dots, z_{n-1})z_n$ for $i < n$ and $z'_n = z_n$, where the $A_i$ are arbitrary local smooth complex-valued functions on $D$. Then the $z'_i$ form a smooth (but not usually holomorphic) complex coordinate system in some small tubular neighbourhood of $D$, and the associated complex coordinate vectors $\{\partial_{z'_i}\}$ have type $(1,0)$ at $z_n = 0$ because $\bar\partial z'_i = O(|z_n|)$. The induced coordinates on $N_D$ are then the same as before because
$\partial_{z'_n} - \partial_{z_n} = \sum_{i < n} A_i \partial_{z_i} \in T^{1,0}D$ at $z_n = 0$.

We observe right away that such a coordinate change has almost no effect at the level of volume forms;
more precisely, there exists a smooth complex $(n-1)$-form $\Upsilon$ such that
\begin{equation}\label{eq:app_vol_form}
dz_1 \wedge \ldots \wedge dz_n = dz'_1 \wedge \ldots \wedge dz'_n +  z_n' \Upsilon \wedge dz_n'.
\end{equation}

The following observation shows how this coordinate freedom can be exploited.

\begin{observation}\label{obs:app_c1_first}
Let $\Phi$ be a smooth map from a neighbourhood of the zero section $D \subset N_D$ to a neighbourhood of $D \subset X$ such that $\Phi(p) = p$ for all points $p \in D$ and such that $d\Phi|_p$ is complex linear with $d\Phi|_p(v) + T^{1,0}_pD = v$ for all points $p \in D$ and vectors $v \in N_{D,p} \subset T_p^{1,0} N_D$. By
changing coordinates on $X$ as above, we can then arrange that
$\Phi^*d{z'_i} = dw_i$ at $w_n = 0$.
\end{observation}

Indeed, always working at the zero section, our assumptions imply that $\Phi^*dz_n = dw_n$, as well as that $\Phi^*dz_i|_{T^{1,0}D} = dw_i|_{T^{1,0}D}$ for all $i < n$. Thus, we simply need to set $A_i = (\Phi^*dz_i)(\partial_{w_n})$.

\begin{example}\label{ex:app_c1}
Fix any background Hermitian metric $g$ on $X$. Let $\Phi$ be the fibrewise $g$-orthogonal projection
$N_D \to (T^{1,0}D) ^\perp$ composed with the $g$-normal exponential map $(T^{1,0}D)^\perp \to X$. Then, for all points $p$ in the zero section, $d\Phi|_p$ is given by $g$-orthogonal projection  composed with the natural isomorphism $T_p^{1,0}(T^{1,0}D)^\perp \to T^{1,0}_pX$. Thus, $\Phi$ satisfies the conditions of Observation \ref{obs:app_c1_first}.
\end{example}

For us, the main use of these constructions lies in the following error estimate.

\begin{lemma}\label{lem:app_c1}
For $\Phi$ and $z'_i$ as in Observation \ref{obs:app_c1_first}, it holds that
\begin{eqnarray}
\label{exp_app}
\Phi^{*}z'_{i}=w_{i}+A_{i,1} w_{n}^{2}+ A_{i,2} w_{n}\bar{w}_n+ A_{i,3}\bar{w}_{n}^{2}.
\end{eqnarray}
Here and in the proof, a letter $A$ with a subscript denotes a generic smooth function.
\end{lemma}

\begin{proof}
The function $\Phi^*z'_i - w_i$ vanishes on the zero section. Thus, \cite[Lemma 2.1]{Milnor} yields that
\begin{equation}\label{hithere}
\Phi^*z'_i - w_i =A_1 w_{n}
+A_2 \bar{w}_{n}.
\end{equation}
By taking the exterior derivative and making use of the fact that $\Phi^* dz'_i = dw_i$ at $w_n = 0$, we derive that
$A_1=A_2=0$ at $w_n = 0$. Thus, using \cite[Lemma 2.1]{Milnor} once again,
$$A_1 =A_3 w_{n}+A_4\bar{w}_{n}\;\,\textrm{and}\;\, A_2=A_5w_{n}+A_6\bar{w}_{n}.$$
Substituting this into \eqref{hithere} immediately yields (\ref{exp_app}).
\end{proof}

\subsection{The test function in Lemma \ref{heythere}}\label{s:test_function} To prove Lemma \ref{heythere}, we were left with checking that
\begin{equation}\label{e:app_c1_basic}
\biggl(\frac{F'(x)}{F(x)}\biggr)^{2}-\frac{F''(x)}{F(x)}
-\frac{F'(x)}{xF(x)}\geq 0\;\,{\rm for\;all}\;x \in (0,a),
\end{equation}
where $F(x) = G(H(x))$ with $G(H) = \frac{H}{1+H}$ and $H(x) = x \exp(\frac{x}{a-x})$. This is {false} for small values of the parameter $a$. However, it is enough to find \emph{some} $a$ that works. Here we will show that any $a \geq 6$ will do. Numerical experiments suggest that the claim is false for $a \leq 1.78$ and true for $a \geq 1.79$.

We begin by computing that
\begin{equation*}
\begin{split}
F'(x) &= \frac{1}{(1+H)^2}\frac{x^2 - ax + a^2}{(x-a)^2}\exp\biggl(\frac{x}{a-x}\biggr),\\
F''(x) &=\frac{1}{(1+H)^2}\frac{1}{(x-a)^4} \biggl(-\frac{2}{x} \frac{H}{1+H} (x^2 - ax + a^2)^2 + (2a^3 - a^2x)\biggr) \exp\biggl(\frac{x}{a-x}\biggr).
\end{split}
\end{equation*}
Thus, the left-hand side of (\ref{e:app_c1_basic}) is given by
\begin{equation}\label{e:app_c1_next}
\frac{1}{H^2(1+H)^2}\frac{1}{(x-a)^4}\biggl[P(x)\exp\left(\frac{x}{a-x}\right) + Q(x)\biggr]\exp\biggl(\frac{2x}{a-x}\biggr),
\end{equation}
where $P(x) = x^5 - ax^4 + 3a^2x^3 - 3a^3x^2 + a^4x$ and $Q(x) = ax^3 - a^3x$.

We only need to prove that the term in square brackets in (\ref{e:app_c1_next}) is nonnegative for all $x \in (0,a)$. We may assume that $a$ is as large as necessary; in fact, sending $x \to 0^+$ shows that we require that $a \geq 1$. Now, $Q(x) < 0$ for all $x \in (0,a)$, so it suffices to prove that
\begin{equation}\label{e:app_c1_red1}
\frac{1}{x}(P(x) + Q(x)) = x^4 - ax^3 + (3a^2 + a)x^2 - 3a^3 x + (a^4 - a^3) \geq 0\;\,{\rm for\;all}\;x \in (0,a).
\end{equation}
To begin with, observe that (\ref{e:app_c1_red1}) would follow from the stronger assertion that
\begin{equation}\label{e:app_c1_red2}
x^4 - ax^3 + 3a^2x^2 - 3a^3 x + (a^4 - a^3) \geq 0\;\,{\rm for\;all}\;x \in \R.
\end{equation}
Next, notice that the quartic polynomial in (\ref{e:app_c1_red2}) has precisely one critical point, $x_0$ say, and that $x_0$ $=$ $\zeta a$, where $\zeta \in (0.5407,0.5408)$ is the unique real solution to $4\zeta^3 - 3\zeta^2  + 6\zeta - 3 = 0$. Consequently, $x_0 \in (0,a)$ is the global minimum of this quartic. The given bounds on $\zeta$ therefore imply that
$$x_0^4 - ax_0^3 + 3a^2x_0^2 - 3a^3 x_0 + a^4 > 0.1819 a^4.$$
Thus (\ref{e:app_c1_red2}), and hence (\ref{e:app_c1_red1}), certainly holds for all $a \geq 5.4976 > (0.1891)^{-1}$.

\subsection{A Gysin sequence for orbifold pairs}\label{s:gysin} The following is needed to prove Proposition \ref{refer}.

\begin{prop}\label{Gysin1}
Let $X$ be a compact complex orbifold and let $D \subset X$ be a suborbifold divisor that contains ${\rm Sing}(X)$. Then there is a long exact sequence of orbifold de Rham cohomology,
\begin{equation}\label{Gysin}
\cdots \to H^{k-2}(D,\R) \to H^{k}(X,\R)\xrightarrow{j^{*}}H^{k}(X\setminus D,\R)\to
H^{k-1}(D,\R)\to\cdots,
\end{equation}
where $j$ denotes the inclusion of $X\setminus D$ into $X$.
\end{prop}

\begin{proof}
Let $U$ be a tubular neighbourhood of $D$. Then we have the usual long exact sequence
\begin{equation*}
\cdots\to H^{k}(X, X\setminus U, \R)\to H^{k}(X,\R)\xrightarrow{j^{*}} H^{k}(X\setminus U,\R)\to H^{k+1}(X, X \setminus U,\R)\to\cdots
\end{equation*}
of orbifold de Rham cohomology, which can obviously be rewritten as
\begin{equation*}
\cdots\to H_{c}^{k}(U,\R)\to H^{k}(X,\R) \xrightarrow{j^{*}} H^{k}(X\setminus D,\R)\to H_{c}^{k+1}(U,\R)\to\cdots.
\end{equation*}
Now define $n = \dim_\C X$ and observe that $H_{c}^{k}(U,\R)=H^{2n-k}(U,\R)
=H^{2n-k}(D,\R)
=H^{k-2}(D,\R)$ by Poincar{\'e} duality and homotopy invariance for orbifold de Rham cohomology \cite{satake}.
\end{proof}

\subsection{Computations for Section \ref{s:setup}}\label{comp4} Let $(X,D)$ denote the pair from Section \ref{s:irregular}.

\subsubsection{Explicit equations and coordinate charts for $X$}\label{charts-mfd} By the computations in \cite[Section 4.2]{dP2-physics}, $X$ is cut out by the following $14$ quadrics in $\P^8$ with coordinates $[z_1:\ldots:z_9]$:
\begin{equation}\label{eqns-mfd}
\begin{split}
&z_{4}^{2}=z_{3}z_{5}\qquad z_{4}^{2}=z_{1}z_{7}+z_{4}z_{9}\qquad z_{4}^{2}=z_{2}z_{6}+z_{4}z_{9}\qquad z_{6}^{2}=z_{3}z_{8}\qquad z_{7}^{2}=z_{5}z_{8}\\
&z_{2}z_{3}=z_{1}z_{4}\qquad z_{4}z_{6}=z_{3}z_{7}\qquad z_{2}z_{4}=z_{1}z_{5}\qquad z_{5}z_{6}=z_{4}z_{7}\qquad z_{3}z_{4}=z_{1}z_{6}+z_{3}z_{9}\\
&z_{4}z_{5}=z_{2}z_{7}+z_{5}z_{9}\qquad z_{6}z_{7}=z_{4}z_{8}\qquad z_{4}z_{6}=z_{1}z_{8}+z_{6}z_{9}\qquad z_{4}z_{7}=z_{2}z_{8}+z_{7}z_{9}.
\end{split}
\end{equation}
It follows directly from this that the affine open sets $X_i = X \cap \{z_i \neq 0\}$ with $i \in \{1,2,3,5,8\}$
cover some neighbourhood of $D = X \cap (z_9 = 0)$; in fact, the union of the $X_i$ is precisely the complement of the point $p_9 = [0:\ldots:0:1]$. We now construct coordinates $(u_i, v_i, w_i)$ on $X_i$ by defining
\begin{equation}\label{indextable}
u_i = \frac{z_{p(i)}}{z_i},\;\, v_i = \frac{z_{q(i)}}{z_i}, \;\, w_i = \frac{z_9}{z_i}, \qquad
\begin{tabular}{c|c c c c c}
$i$ & 1 & 2 & 3 & 5 & 8 \\\hline
$p(i)$ & $2$ & $1$ & $1$ & $2$ & $6$  \\
$q(i)$ & $3$ & $5$ & $6$ & $7$ & $7$  \\
\end{tabular}.
\end{equation}
To show that these are indeed coordinates (and hence that $X_i$ is isomorphic to $\C^3$), we simply write down explicit formulas for the inverse maps from $\C^3$ into $X$ as follows:
\begin{equation*}
\begin{split}
&(u_{1},v_{1},w_{1}) \mapsto [1:u_{1}:v_{1}:u_{1}v_{1}:u_{1}^{2}v_{1}:v_{1}(u_{1}v_{1}-w_{1}):u_{1}v_{1}(u_{1}v_{1}-w_{1}):v_{1}
(u_{1}v_{1}-w_{1})^{2}:w_{1}],\\
&(u_2,v_2,w_2) \mapsto [u_{2}:1:u_{2}^{2}v_{2}:u_{2}v_{2}:v_{2}:u_{2}v_{2}(u_{2}v_{2}-w_{2}):v_{2}
(u_{2}v_{2}-w_{2}):v_{2}(u_{2}v_{2}-w_{2})^{2}:w_{2}],\\
&(u_3,v_3,w_3)\mapsto [u_{3}:u_{3}(u_{3}v_{3}+w_{3}):1:u_{3}v_{3}+w_{3}:(u_{3}v_{3}+w_{3})^{2}:v_{3}:v_{3}(u_{3}v_{3}+w_{3}):v_{3}
^{2}:w_{3}],\\
&(u_5, v_5, w_5) \mapsto [u_{5}(u_{5}v_{5}+w_{5}):u_{5}:(u_{5}v_{5}+w_{5})^{2}:u_{5}v_{5}+w_{5}:1:v_{5}(u_{5}v_{5}+w_{5}):v_{5}:v_{5}
^{2}:w_{5}],\\
&(u_8, v_8, w_8) \mapsto [u_{8}(u_{8}v_{8}-w_{8}):v_{8}(u_{8}v_{8}-w_{8}):u_{8}^{2}:u_{8}v_{8}:v_{8}^{2}:u_{8}:v_{8}:1:w_{8}].
\end{split}
\end{equation*}

As mentioned in Section \ref{4:charts:X}, using these formulas, one can readily check  that the locally defined meromorphic $(3,0)$-forms $\sigma(i)  du_i \wedge dv_i \wedge dw_i/w_i^2$, where $\sigma(i) = +1$ for $i = 1,3,8$ and $\sigma(i) = -1$ for $i = 2,5$, patch up as a global meromorphic $(3,0)$-form $\Omega$ on $X \setminus \{p_9\}$ with double poles along $D$ and no other poles or zeros; of course, $\Omega$ then extends from $X \setminus \{p_9\}$ to the whole of $X$.

For the sake of completeness, let us also quickly explain how to construct a coordinate chart near the missing point $p_9$, thus confirming that $X$ is indeed smooth everywhere. However, unlike the $X_i$,
this chart will only be isomorphic to a polydisk. The coordinates are given by
$(u,v,w) = (\frac{z_1}{z_9},\frac{z_2}{z_9},\frac{z_8}{z_9})$. To invert this map locally near $p_9$, let $\sigma = f(\tau)$ denote the local inverse to $\tau = \sigma(\sigma-1)^2$ near  $\sigma = 0$. We set $z_9 = 1$ and attempt to determine $z_3,z_4,z_5,z_6,z_7$ from $(u,v,w) = (z_1,z_2,z_8)$. By the 2nd and 14th equations of (\ref{eqns-mfd}), $z_4 = f(uvw)$ provided that $|z_4| \ll 1$. Given this, we can compute $z_6, z_7$ and $z_3, z_5$ from equations $\# 13, 14$ and $\# 10,11$ of (\ref{eqns-mfd}), respectively.

\subsubsection{Explicit equations and coordinate charts for the normal bundle to $D$ in $X$}\label{charts-normal} Again from \cite{dP2-physics},
the following cone in $\C^8$ with its canonical $\C^*$-action $(t,Z) \mapsto tZ$ is $\C^*$-equivariantly isomorphic to $K_D^\times$, the canonical bundle of $D$ with its zero section blown down:
 \begin{equation}\label{eqns-cone}
\begin{matrix}
&\quad Z_{4}^{2}=Z_{3}Z_{5} & & \quad Z_{4}^{2}=Z_{1}Z_{7} & & \quad Z_{4}^{2}=Z_{2}Z_{6} & & \quad Z_{6}^{2}=Z_{3}Z_{8} & & \quad Z_{7}^{2}=Z_{5}Z_{8}\\
&Z_{2}Z_{3}=Z_{1}Z_{4} & & Z_{4}Z_{6}=Z_{3}Z_{7} & & Z_{2}Z_{4}=Z_{1}Z_{5} & & Z_{5}Z_{6}=Z_{4}Z_{7} & & Z_{3}Z_{4}=Z_{1}Z_{6}\\
&Z_{4}Z_{5}=Z_{2}Z_{7} & & Z_{6}Z_{7}=Z_{4}Z_{8} & & Z_{4}Z_{6}=Z_{1}Z_{8} & & \;Z_{4}Z_{7}=Z_{2}Z_{8}. & &
\end{matrix}
\end{equation}
Observe that (\ref{eqns-cone}) can be derived from (\ref{eqns-mfd}) by dehomogenising $z_i = Z_i z_9$ and dropping all linear terms from the resulting system of affine quadrics. We now pass back to the projective completion of (\ref{eqns-cone}) in $\P^8$ since we are mainly interested in a neighbourhood of the divisor at infinity.

The completion of (\ref{eqns-cone}) in $\P^8$ is isomorphic to the compactification of the total space of $N_{D}$, obtained by adding $p_9 = [0:\ldots:0:1]$. Unlike in Section \ref{charts-mfd}, this point is now singular. However, we still have a covering by coordinate charts
$(U_i, V_i, W_i)$ off of $p_9$, where $U_i = {z_{p(i)}}/{z_9}$, $V_i = z_{q(i)}/z_9$, $W_i = z_9/z_i$, and $p(i), q(i)$ are as in (\ref{indextable}). The inverse maps from $\C^3$ into $N_{D}$ are given by
\begin{equation*}
\begin{split}
&(U_{1},V_{1},W_{1}) \mapsto [1:U_1:V_1:U_{1}V_{1}:U_{1}^{2}V_{1}:U_1V_1^2:U_1^2V_1^2:U_1^2V_1^3:W_1],\\
&(U_2,V_2,W_2) \mapsto [U_{2}:1:U_{2}^{2}V_{2}:U_{2}V_{2}:V_{2}:U_{2}^2V_{2}^2:U_2V_{2}^2:U_2^2V_2^3:W_2],\\
&(U_3,V_3,W_3)\mapsto [U_{3}:U_{3}^2V_3:1:U_{3}V_{3}:U_{3}^2V_{3}^2:V_{3}:U_{3}V_3^2:V_3^2:W_3],\\
&(U_5, V_5, W_5) \mapsto [U_{5}^2V_{5}: U_{5}: U_{5}^2V_{5}^2: U_{5}V_{5}: 1: U_5V_{5}^2: V_{5}: V_{5}
^{2}:W_5],\\
&(U_8, V_8, W_8) \mapsto [U_{8}^2V_8: U_{8}V_{8}^2: U_{8}^{2}: U_{8}V_{8}: V_{8}^{2}: U_{8}: V_{8}: 1:W_8].
\end{split}
\end{equation*}

As in Section \ref{charts-mfd}, it follows from these formulas that the locally defined meromorphic $3$-forms  $\sigma(i) {dU_{i}\wedge dV_{i}\wedge dW_{i}}/{W_{i}^{2}}$ patch up as a global meromorphic $3$-form $\Omega_0$ on $N_{D}$, homogeneous under the canonical $\C^*$-action on $N_{D}$ and with double poles
along the zero section. \newpage

\bibliographystyle{amsplain}
\bibliography{ref}

\end{document}